\documentclass[11pt]{article}

\usepackage{amsmath, amssymb, amscd}
\usepackage[matrix,arrow]{xy}

\def\eqref#1{(\ref{#1})}
\newcommand{\goth}{\mathfrak}

\newcommand{\arrow}{{\:\longrightarrow\:}}
\newcommand{\Z}{{\Bbb Z}}
\newcommand{\C}{{\Bbb C}}
\newcommand{\R}{{\Bbb R}}
\newcommand{\Q}{{\Bbb Q}}

\newcommand{\6}{\partial}
\def\1{\sqrt{-1}\:}
\newcommand{\restrict}[1]{{\left|_{{\phantom{|}\!\!}_{#1}}\right.}}
\newcommand{\cntrct}                % contraction with a vector field
{\hspace{2pt}\raisebox{1pt}{\text{$\lrcorner$}}\hspace{2pt}}

\renewcommand{\tilde}{\widetilde}
\renewcommand{\bar}{\overline}
\renewcommand{\phi}{\varphi}
\renewcommand{\epsilon}{\varepsilon}
\renewcommand{\geq}{\geqslant}
\renewcommand{\leq}{\leqslant}

\newcommand{\End}{\operatorname{End}}

\newcommand{\Id}{\operatorname{Id}}

\newcommand{\Aut}{\operatorname{Aut}}

\newcommand{\codim}{\operatorname{codim}}

\newcommand{\comment}[1]{{}}

\def\blacksquare{\hbox{\vrule width 4pt height 4pt depth 0pt}}
\def\endproof{\blacksquare}

\makeatletter

\@ifundefined{Bbb}
     {\newcommand{\Bbb}[1]{{\mathbb #1}}}%
{}%     {\edef\Bbb#1{{\Bbb #1}}}

\newcommand{\ps@verbit}{%
  \renewcommand{\@oddhead}{%
          \scriptsize
          {Vanishing theorems for hyperk\"ahler manifolds}
          \hfil\tiny {M. Verbitsky,  \ \ \ \ April 9, 2006}}
  \renewcommand{\@evenhead}{\@oddhead}
  \renewcommand{\@oddfoot}{\hfil\thepage\hfil}
  \renewcommand{\@evenfoot}{\@oddfoot}}
 
\newcounter{Mycounter}[section]
\newcounter{lemma}[section]
\renewcommand{\thelemma}{{Lemma \thesection.\arabic{lemma}}}
\newcommand{\lemma}{%
     \setcounter{lemma}{\value{Mycounter}}
     \refstepcounter{lemma}
     \stepcounter{Mycounter}
     {\bf \thelemma:\ }}

\newcounter{claim}[section]
\renewcommand{\theclaim}{{Claim \thesection.\arabic{claim}}}
\newcommand{\claim}{%
     \setcounter{claim}{\value{Mycounter}}
     \refstepcounter{claim}
     \stepcounter{Mycounter}
     {\bf \theclaim:\ }}

\newcounter{sublemma}[section]
\newcounter{corollary}[section]
\renewcommand{\thecorollary}{{Corollary \thesection.\arabic{corollary}}}
\newcommand{\corollary}{%
     \setcounter{corollary}{\value{Mycounter}}
     \refstepcounter{corollary}
     \stepcounter{Mycounter}
     {\bf \thecorollary:\ }}

\newcounter{theorem}[section]
\renewcommand{\thetheorem}{{Theorem \thesection.\arabic{theorem}}}
\newcommand{\theorem}{%
     \setcounter{theorem}{\value{Mycounter}}
     \refstepcounter{theorem}
     \stepcounter{Mycounter}
     {\bf \thetheorem:\ }}

\newcounter{conjecture}[section]
\newcounter{proposition}[section]
\renewcommand{\theproposition}
       {{Proposition \thesection.\arabic{proposition}}}
\newcommand{\proposition}{%
     \setcounter{proposition}{\value{Mycounter}}
     \refstepcounter{proposition}
     \stepcounter{Mycounter}
     {\bf \theproposition:\ }}

\newcounter{definition}[section]
\renewcommand{\thedefinition}
       {{Definition~\thesection.\arabic{definition}}}
\newcommand{\definition}{%
     \setcounter{definition}{\value{Mycounter}}
     \refstepcounter{definition}
     \stepcounter{Mycounter}
     {\bf \thedefinition:\ }}

\newcounter{example}[section]
\newcounter{remark}[section]
\renewcommand{\theremark}{{Remark \thesection.\arabic{remark}}}
\newcommand{\remark}{%
     \setcounter{remark}{\value{Mycounter}}
     \refstepcounter{remark}
     \stepcounter{Mycounter}
     {\bf \theremark:\ }}

\newcounter{problem}[section]
\newcounter{question}[section]
\@addtoreset{equation}{section}
\@addtoreset{footnote}{section}
\makeatother

\begin{document}

%%%%%%%%%%%%%%%%%%%%%%%%%%%%%%%%%%%%%%%%%%%%%%%%%%%%%%%%%%%%
\begin{center}
{\Large\bf
Quaternionic Dolbeault complex and vanishing \\[3mm]theorems
on hyperk\"ahler manifolds
}
%%%%%%%%%%%%%%%%%%%%%%%%%%%%%%%%%%%%%%%%%%%%%%%%%%%%%%%%%%%%
\\[4mm]
Misha Verbitsky,\footnote{The author is 
supported by  EPSRC grant  GR/R77773/01}
\\[4mm]
{\tt verbit@maths.gla.ac.uk, verbit@mccme.ru}
\end{center}

%%%%%%%%%%%%%%%%%%%%%%%%%%%%%%%%%%%%%%%%%%%%%%%%
{\small 
\hspace{0.15\linewidth}
\begin{minipage}[t]{0.7\linewidth}
{\bf Abstract} \\
Let $(M,I,J,K)$ be a compact hyperk\"ahler manifold,
$\dim_{\Bbb H}M=n$, and $L$ a non-trivial 
holomorphic line bundle on $(M,I)$. Using the 
quaternionic Dolbeault complex, we prove the
following vanishing theorem for holomorphic 
cohomology of $L$.  If $c_1(L)$ 
lies in the closure $\hat K$ of the dual 
K\"ahler cone, then $H^i(L)=0$ for $i>n$. 
If $c_1(L)$ lies in the opposite cone $-\hat K$, 
then $H^i(L)=0$ for $i<n$. Finally, if $c_1(L)$
is neither in $\hat K$ nor in $-\hat K$, then
$H^i(L)=0$ for $i\neq n$. 
\end{minipage}
}
%%%%%%%%%%%%%%%%%%%%%%%%%%%%%%%%%%%%%%%%%%%%%%%%

{
\small
\tableofcontents
}
%%%%%%%%%%%%%%%%%%%%%%%%%%%%%%%%%%%%%%%%%%%%%%%%

\section{Introduction}
\label{_Intro_Section_}

%%%%%%%%%%%%%%%%%%%%%%%%%%%%%%%%%%%%%%%%%%%%%%%%

%%%%%%%%%%%%%%%%%%%%%%%%%%%%%%%%%%%%%%%%%%%%%%%%
\subsection{Hypercomplex and hyperk\"ahler manifolds}
\label{_HK_intro_Subsection_}
%%%%%%%%%%%%%%%%%%%%%%%%%%%%%%%%%%%%%%%%%%%%%%%%

%%%%%%%%%%%%%%%%%%%%%%%%%%%%%%%%%%%%%%%%%%%%%%%%
\definition
Let $M$ be a manifold, and $I,J,K\in \End(TM)$
endomorphisms of the tangent bundle satisfying the
quaternionic relation
\[
I^2=J^2=K^2=IJK=-\Id_{TM}.
\]
The manifold $(M,I,J,K)$ is called {\bf hypercomplex}
if the almost complex structures $I$, $J$, $K$
are integrable. If, in addition, $M$
is equipped with a Riemannian metric $g$ which
is K\"ahler with respect to $I,J,K$, the
manifold $(M,I,J,K,g)$
is called {\bf hyperk\"ahler}. 

Consider the K\"ahler forms $\omega_I, \omega_J, \omega_K$
on $M$:
\[
\omega_I(\cdot, \cdot):= g(\cdot, I\cdot), \ \
\omega_J(\cdot, \cdot):= g(\cdot, J\cdot), \ \
\omega_K(\cdot, \cdot):= g(\cdot, K\cdot).
\]
An elementary linear-algebraic calculation implies
that the 2-form $\Omega:=\omega_J+\1\omega_K$ is of Hodge type $(2,0)$
on $(M,I)$. This form is clearly closed and
non-degenerate, hence it is a holomorphic symplectic form.

In algebraic geometry, the word ``hyperk\"ahler''
is essentially synonymous with ``holomorphically
symplectic'', due to the following theorem, which is
implied by Yau's solution of Calabi conjecture.

\hfill

%%%%%%%%%%%%%%%%%%%%%%%%%%%%%%%%%%%%%%%%%%%%%%%%
\theorem\label{_Calabi-Yau_Theorem_}
Let $(M,I)$ be a compact, K\"ahler, holomorphically
symplectic manifold. Then there exists a unique
hyperk\"ahler metric on $(M,I)$ with the same K\"ahler class.

\hfill

{\bf Proof:} See \cite{_Yau:Calabi-Yau_}, \cite{_Besse:Einst_Manifo_}.
\endproof

\hfill

%%%%%%%%%%%%%%%%%%%%%%%%%%%%%%%%%%%%%%%%%%%%%%%%
\remark
The hyperk\"ahler metric is unique, but there could
be several hyperk\"ahler structures compatible with
a given hyperk\"ahler metric on $(M,I)$.

%%%%%%%%%%%%%%%%%%%%%%%%%%%%%%%%%%%%%%%%%%%%%%%%%%%%%%%%%%%%%%
\subsection{Bogomolov's decomposition theorem}
%%%%%%%%%%%%%%%%%%%%%%%%%%%%%%%%%%%%%%%%%%%%%%%%%%%%%%%%%%%%%%

The modern approach to Bogomolov's decomposition
is based on Calabi-Yau theorem
(\ref{_Calabi-Yau_Theorem_}), Berger's classification
of irreducible holonomy and de Rham's splitting theorem
for holonomy reduction (\cite{_Beauville_},
\cite{_Besse:Einst_Manifo_}). It is worth mention
that the original proof of decomposition theorem
(due to F. Bogomolov, \cite{_Bogomolov:decompo_}) 
was algebraic.

\hfill

%%%%%%%%%%%%%%%%%%%%%%%%%%%%%%%%%%%%%%%%%%%%%%%%%%%%%%%%
\theorem\label{_Bogomo_splitting_} 
Let $(M,I,J,K)$ be a compact hyperk\"ahler manifold.
Then there exists a finite covering $\tilde M\arrow M$,
such that $\tilde M$ is decomposed, as a hyperk\"ahler
manifold, into a product
\[
\tilde M = M_1\times M_2 \times \dots M_n \times T,
\]
where $(M_i,I,J,K)$ satisfy $H^1(M_i)=0$, $H^{2,0}(M_i,I)=\C$,
and $T$ is a hyperk\"ahler torus. Moreover,
$M_i$ are uniquely determined by $M$ and simply
connected, and $T$ is unique up to isogeny.

\hfill

{\bf Proof:} See \cite{_Beauville_},
\cite{_Besse:Einst_Manifo_}.\endproof

\hfill

%%%%%%%%%%%%%%%%%%%%%%%%%%%%%%%%%%%%%%%%%%%%%%%%%%%%%%%%
\definition
Let $(M,I,J,K)$ be a compact hyperk\"ahler manifold which satisfies 
$H^1(M)=0$, $H^{2,0}(M,I)=\C$. Then $M$ is called
{\bf a irreducible hyperk\"ahler manifold}.

\hfill

%%%%%%%%%%%%%%%%%%%%%%%%%%%%%%%%%%%%%%%%%%%%%%%%%%%%%%%%
\remark
Notice that \ref{_Bogomo_splitting_} implies that 
irreducible hyperk\"ahler manifolds are simly connected.
In particular, they do not admit a further decomposition.
This explains the term ``irreducible''.

%%%%%%%%%%%%%%%%%%%%%%%%%%%%%%%%%%%%%%%%%%%%%%%%
\subsection{Vanishing theorems on hyperk\"ahler manifolds}
%%%%%%%%%%%%%%%%%%%%%%%%%%%%%%%%%%%%%%%%%%%%%%%%

Using the argument which essentially belongs to the 
theory of hypercomplex manifolds, we are able to prove the
following algebro-geometric statements.

\hfill

%%%%%%%%%%%%%%%%%%%%%%%%%%%%%%%%%%%%%%%%%%%%%%%%%%%%%%%%
\theorem\label{_main_vanishing_intro_Theorem_}
Let $(M,I,J,K)$ be a compact, irreducible hyperk\"ahler
manifold, and $L$ a holomorphic line bundle on $(M,I)$ with
$c_1(L)\neq 0$. Denote by 
\[ \bar{\cal K}\check{\;}\subset 
    H^{1,1}(M,I)\cap H^2(M,\R)
\]
the closure of the dual K\"ahler cone of $(M,I)$
(Subsection \ref{_dual_Kahler_Subsection_}), 
and let $-\bar{\cal K}\check{\;}$
be the opposite cone. Then one of the following holds.
\begin{description}
\item[(i)] $c_1(L)\in \bar{\cal K}\check{\;}$; then
$H^i(L) =0$ for $i> \frac{\dim_\C M}{2}$.
\item[(ii)] $c_1(L)\in -\bar{\cal K}\check{\;}$; then
$H^i(L) =0$ for $i< \frac{\dim_\C M}{2}$.
\item[(iii)] $c_1(L)$ does not lie in
$-\bar{\cal K}\check{\;}\cup \bar{\cal K}\check{\;}$; then
$H^i(L) =0$ for $i\neq \frac{\dim_\C M}{2}$.
\end{description}
{\bf Proof:} See
\ref{_main_vanishing_after_cones_Theorem_}.
\endproof

\hfill

%%%%%%%%%%%%%%%%%%%%%%%%%%%%%%%%%%%%%%%%%%%%%%%%%%%%%%%%%%%%%%
\theorem\label{_vanishing_arbi_rank_intro_Theorem_}
Let $(M,I,J,K)$ be a compact, irreducible hyperk\"ahler
manifold, $L$ a holomorphic line bundle on $(M,I)$ with
$c_1(L)\neq 0$, and $B$ an arbitrary holomorphic vector
bundle on $(M,I)$. Then there exists a sufficiently
big number $N_0$, such that for any integer $N>N_0$
one of the following holds.
\begin{description}
\item[(i)] $c_1(L)\in \bar{\cal K}\check{\;}$; then
$H^i(L^N\otimes B) =0$ for $i> \frac{\dim_\C M}{2}$.
\item[(ii)] $c_1(L)\in -\bar{\cal K}\check{\;}$; then
$H^i(L^N\otimes B) =0$ for $i< \frac{\dim_\C M}{2}$.
\item[(iii)] $c_1(L)$ does not lie in
$-\bar{\cal K}\check{\;}\cup \bar{\cal K}\check{\;}$; then
$H^i(L^N\otimes B) =0$ for $i\neq \frac{\dim_\C M}{2}$.
\end{description}
{\bf Proof:} This is \ref{_vanishing_arb_rank_Theorem_}.
\endproof

\hfill

The vanishing theorems have many interesting
geometrical consequences. As an example, we give
the following theorem (Section \ref{_vanishing_and_nef_Section_}).

\hfill

%%%%%%%%%%%%%%%%%%%%%%%%%%%%%%%%%%%%%%%%%%%%%%%%
\theorem\label{_Surjecti_comple_inte_Theorem_}
Let $(M, I, J, K)$ be a irreducible hyperk\"ahler manifold,
and $X\subset (M,I)$ a subvariety of dimension
$\dim_\C X>\frac 1 2 \dim_\C M$. Assume that
$X$ is a complete intersection of ample divisors.
Consider a holomorphic line bundle $L$ on
$(M,I)$ with $c_1(L)$ nef (that is, $c_1(L)$
lies in the closure of the K\"ahler cone of $(M,I)$) and
$q(c_1(L), c_1(L))=0$, where $q$ is the 
Bogomolov-Beauville-Fujiki bilinear form 
(\ref{_BBF_Definition_}). Then the natural
restriction map is surjective on holomorphic sections:
\[
H^0(L^N) \arrow H^0(L^N\restrict X)\arrow 0.
\]
for a sufficiently big power of $L$.

\hfill

{\bf Proof:} See \ref{_comple_interse_surj_Theorem_}. \endproof

%%%%%%%%%%%%%%%%%%%%%%%%%%%%%%%%%%%%%%%%%%%%%%%%%%%%%%%%%%%%%%
\subsection{Quaternionic Dolbeault complex and vanishing}
%%%%%%%%%%%%%%%%%%%%%%%%%%%%%%%%%%%%%%%%%%%%%%%%%%%%%%%%%%%%%%

In this Subsection we give a brief introduction
to quaternionic Dolbeault complex. We sketch
how one can use the quaternionic Dolbeault complex
to deduce the vanishing theorems for cohomology.
Further on in this paper, this theme is developed
in a more detailed way.

Let $M$ be a hypercomplex manifold. 
There is a natural action of $SU(2)$ on $\Lambda^1(M)$
(we identify $SU(2)$ with the group of unitary
quaternions). This action is extended to $\Lambda^*(M)$
by multiplicativity.

This $SU(2)$-action plays the same role 
in hypercomplex and hyperk\"ahler
geometry as the usual Hodge decomposition
 in complex geometry.

Let 
\[
\Lambda^i(M) = \bigoplus \Lambda^i_k(M)
\]
be a weight decomposition of the space of $i$-forms,
with $\Lambda^i_k(M)$ an $SU(2)$-representation of weight
$k$ (see Subsection \ref{_Weights_Subsection_}).
It is easy to check that 
\[
V^*:= \bigoplus_{i>k}\Lambda^i_k(M)
\]
is a differential ideal in the de Rham algebra 
$\Lambda^*(M)$, that is, an ideal which satisfies $dV^*\subset V^*$
(Subsection \ref{_qD_Subsection_}).
Therefore, the quotient $\Lambda^*(M)/V^*$ 
is a differential graded algebra, denoted as
$(\Lambda^*_+(M), d_+)$. This algebra is
called {\bf the quaternionic Dolbeault complex} 
(\ref{_qD_Definition_}). We approach 
$(\Lambda^*_+(M), d_+)$ from the same point of view
as one approaches the usual Dolbeault complex
in algebraic geometry. There is a Hodge decomposition
(Subsection \ref{_Hodge_on_qD_Subsection_}),
and a Lefschetz-type ${\goth{sl}}(2)$-action
(\ref{_SL(2)_Proposition_}). The analogue
of Kodaira-Nakano formula is written in  
\eqref{_Laplacians_K-N_ide_Equation_}:
\begin{equation}\label{_KN_Intro_Equation_}
\Delta_{\bar\6}= \Delta_{\bar\6_J}+ [ \Lambda_{\bar\Omega},\Theta_+],
\end{equation}
where $\Delta_{\bar\6}=\bar\6\bar\6^*+\bar\6^*\bar\6$
is the usual Laplacian on $(0,p)$-forms with coefficients
in a holomorphic vector bundle $B$ on $(M,I)$,
$\Theta_+$ the $\Lambda^2_+(M)\otimes \End(B)$-part
of the curvature of $B$, and $\Delta_{\bar\6_J}$
a positive self-dual operator. When the commutator
$[ \Lambda_{\bar\Omega},\Theta_+]$ is positive, this
leads to the vanishing theorems
\ref{_main_vanishing_intro_Theorem_} and 
\ref{_vanishing_arbi_rank_intro_Theorem_}, which are
deduced from \eqref{_KN_Intro_Equation_} in the
same way as Kodaira-Nakano vanishing is deduced
from the Kodaira-Nakano identity.

If the bundle $B$ is a line bundle,
we can choose its metric in such a way that
its curvature 2-form $\Theta_B$ is harmonic
(\cite{_Griffi_Harri_}).
Consider the weight decomposition
\[
\Theta_B=\Theta_+  + \Theta_0,
\] 
where $\Theta_0$ is $SU(2)$-invariant.
Then $\Theta_+$ is harmonic (see Subsection
\ref{_SU(2)_on_coho_Subsection_}).
From \eqref{_primiti_decompo_H^2_Equation_}, 
it follows that $\Theta_+=\lambda\omega_I$,
where $\lambda$ is a real constant, and $\omega_I$
is the K\"ahler form of $(M,I)$.
Then 
\begin{equation}\label{_commu_via_H_intro_Equation_}
[ \Lambda_{\bar\Omega},\Theta_+]= \lambda H_{\bar\Omega},
\end{equation}
 where $H$ is a scalar operator mapping a $(0,p)$-form $\eta$
into $(n -p)\eta$, where $n=\dim_{\Bbb H} M$
(see \eqref{_H_as_commu_and_Lapla_Equation_}). 
For $\lambda>0$, \eqref{_commu_via_H_intro_Equation_}
is positive when $p <n$, and for $\lambda<0$,
\eqref{_commu_via_H_intro_Equation_}
is positive when $p >n$. The vanishing
of holomorphic cohomology
(for $p>n$ in the first case, and for
$p<n$ in the second case) is a consequence.

%%%%%%%%%%%%%%%%%%%%%%%%%%%%%%%%%%%%%%%%%%%%%%%%

\section{Quaternionic Dolbeault complex}
\label{_QD_Section_}

%%%%%%%%%%%%%%%%%%%%%%%%%%%%%%%%%%%%%%%%%%%%%%%%

%%%%%%%%%%%%%%%%%%%%%%%%%%%%%%%%%%%%%%%%%%%%%%%%
\subsection{Weights of $SU(2)$-representations}
\label{_Weights_Subsection_}
%%%%%%%%%%%%%%%%%%%%%%%%%%%%%%%%%%%%%%%%%%%%%%%%

It is well-known that any irreducible representations
of $SU(2)$ over $\C$ can be obtained as a symmetric power
$S^i(V_2)$, where $V_1$ is a fundamental 2-dimensional
representation. We say that a representation $W$ 
{\bf has weight $i$} if it is isomorphic to $S^i(V_1)$.
A representation is said to be {\bf pure of weight $i$}
if all its irreducible components have weight $i$.
If all irreducible components of a representation $W_1$
have weight $\leq i$, we say that $W_1$ {\bf is a
  representation of weight $\leq i$}.
In a similar fashion one defines representations
of weight $\geq i$.

\hfill

%%%%%%%%%%%%%%%%%%%%%%%%%%%%%%%%%%%%%%%%%%%%%%%%
\remark\label{_weight_multi_Remark_}
The Clebsch-Gordan formula (see \cite{_Humphreys_})
claims that the weight is {\em multiplicative}, 
in the following sense: if $i\leq j$, then
\[
V_i\otimes V_j = \bigoplus_{k=0}^i V_{i+j-2k},
\]
where $V_i=S^i(V_1)$ denotes the irreducible
representation of weight $i$.

\hfill

A subspace $W\subset W_1$ is {\bf pure of weight $i$}
if the $SU(2)$-representation $W'\subset W_1$ generated
by $W$ is pure of weight $i$.

%%%%%%%%%%%%%%%%%%%%%%%%%%%%%%%%%%%%%%%%%%%%%%%%
\subsection{Quaternionic Dolbeault complex: a definition}
\label{_qD_Subsection_}
%%%%%%%%%%%%%%%%%%%%%%%%%%%%%%%%%%%%%%%%%%%%%%%%

Let $M$ be a hypercomplex (e.g. a hyperk\"ahler) manifold,
$\dim_{\Bbb H}M=n$.
There is a natural multiplicative action of $SU(2)\subset
{\Bbb H}^*$ on $\Lambda^*(M)$, associated with the
hypercomplex structure. 

\hfill

%%%%%%%%%%%%%%%%%%%%%%%%%%%%%%%%%%%%%%%%%%%%%%%%
\remark\label{_weights_on_forms_Remark_}
The space $\Lambda^*(M)$
is an infinite-dimensional representation of $SU(2)$,
however, all its irreducible components
are finite-dimensional. Therefore it makes
sense to speak of {\em weight} of $\Lambda^*(M)$
and its sub-\-rep\-re\-sen\-ta\-tions. Clearly, $\Lambda^1(M)$
has weight 1. From Clebsch-Gordan formula
(\ref{_weight_multi_Remark_}), it follows that 
 $\Lambda^i(M)$ is an $SU(2)$-representation
of weight $\leq i$. Using the Hodge $*$-isomorphism 
$\Lambda^i(M)\cong \Lambda^{4n-i}(M)$, we find that
for $i> 2n$, $\Lambda^i(M)$ is a representation
of weight $\leq 2n-i$.

\hfill

Let $V^i\subset \Lambda^i(M)$ be a maximal
$SU(2)$-invariant subspace of weight $<i$.
The space $V^i$ is well defined, because
it is a sum of all irreducible representations
$W\subset \Lambda^i(M)$ of weight $<i$.
Since the weight is multiplicative
(\ref{_weight_multi_Remark_}), $V^*= \bigoplus_i V^i$
is an ideal in $\Lambda^*(M)$. We also have
$V^i = \Lambda^i(M)$ for $i> 2n$
(\ref{_weights_on_forms_Remark_}).

It is easy to see that the de Rham differential
$d$ increases the weight by 1 at most. Therefore,
$dV^i\subset V^{i+1}$, and $V^*\subset \Lambda^*(M)$
is a differential ideal in the de Rham DG-algebra
$(\Lambda^*(M), d)$.

\hfill

%%%%%%%%%%%%%%%%%%%%%%%%%%%%%%%%%%%%%%%%%%%%%%%%
\definition\label{_qD_Definition_}
Denote by $(\Lambda^*_+(M), d_+)$ the quotient algebra
$\Lambda^*(M)/V^*$
It is called {\bf the quaternionic Dolbeault algebra of
  $M$}, or {\bf the quaternionic Dolbeault complex} 
(qD-algebra or qD-complex for short).

The space $\Lambda^i_+(M)$ can be identified with the
maximal subspace of $\Lambda^i(M)$ of weight $i$,
that is, a sum of all irreducible sub-representations of weight $i$.
This way, $\Lambda^i_+(M)$ can be considered as a subspace
in $\Lambda^i(M)$; however, this subspace is not preserved
by the multiplicative structure and the differential.

\hfill

%%%%%%%%%%%%%%%%%%%%%%%%%%%%%%%%%%%%%%%%%%%%%%%%
\remark 
The complex $(\Lambda^*_+(M), d_+)$ 
was constructed much earlier by Salamon,
in a different (and much more general) situation,
and much studied since then
(\cite{_Salamon_}, \cite{_Capria-Salamon_},
\cite{_Baston_}, \cite{_Leung_}).

%%%%%%%%%%%%%%%%%%%%%%%%%%%%%%%%%%%%%%%%%%%%%%%%
\subsection{The Hodge decomposition on the quaternionic
  Dolbeault complex}.
\label{_Hodge_on_qD_Subsection_}
%%%%%%%%%%%%%%%%%%%%%%%%%%%%%%%%%%%%%%%%%%%%%%%%

Let $(M,I,J,K)$ be a hypercomplex 
manifold, and $L$ a complex structure
 induced by the quaternionic action, say,
$I$, $J$ or $K$.
 Consider the $U(1)$-action
on $\Lambda^1(M)$ provided by 
$\phi\stackrel{\rho_L} \arrow \cos \phi\Id + \sin\phi \cdot L$.
We extend this action to a multiplicative action on
$\Lambda^*(M)$. Clearly, for a $(p,q)$-form 
$\eta\in \Lambda^{p,q}(M,L)$, we have
\begin{equation}\label{_Hodge_weights_Equation_}
   \rho_L(\phi)\eta = e^{\1(p-q)\phi}\eta.
\end{equation}

\hfill

%%%%%%%%%%%%%%%%%%%%%%%%%%%%%%%%%%%%%%%%%%%%%%%%
\lemma \label{_su(2)_action_explici_Lemma_}
Let $(M,I,J,K)$ be a hypercomplex manifold
and \[ \rho_I, \rho_J, \rho_K\] the 
homomorphisms \[ U(1) \arrow \Aut(\Lambda^*(M))\] constructied above. 
Then $\rho_I, \rho_J, \rho_K$ generate 
the Lie group action \[ SU(2)\subset \Aut(\Lambda^*(M))\]
associated with the hypercomplex structure.

\hfill

{\bf Proof:} \ref{_su(2)_action_explici_Lemma_} is clear.
Indeed, the action of $SU(2)$, and $\rho_I, \rho_J,
\rho_K$  are defined on $\Lambda^*(M)$ by
multiplicativity, hence it suffices to check that
$\rho_I, \rho_J, \rho_K$ generate the standard action of
$SU(2)$ on $\Lambda^1(M)$. On $\Lambda^1(M)$,
$\rho_I, \rho_J, \rho_K$ act as quaternions
$\cos \phi + \sin\phi \cdot I$, $\cos \phi + \sin\phi
\cdot J$, $\cos \phi + \sin\phi \cdot K$, and
they generate the group of unitary quaternions.
\endproof

\hfill

From \ref{_su(2)_action_explici_Lemma_},
it is clear that $\rho_L$ preserves 
components of weight $i$. We obtain that
$V^*$ is preserved by $\rho_L$, hence $\rho_L$
acts on $\Lambda^*_+(M)$. Then,
\eqref{_Hodge_weights_Equation_}
gives a Hodge decomposition on $\Lambda^*_+(M)$:
\[
\Lambda^i_+(M) = \bigoplus_{p+q=i}\Lambda^{p,q}_{+,L}(M).
\]

The following result is implied immediately by the
standard calculations from the theory of
$SU(2)$-representations.

\hfill

%%%%%%%%%%%%%%%%%%%%%%%%%%%%%%%%%%%%%%%%%%%%%%%%
\proposition \label{_qD_decompo_expli_Proposition_}
Let $(M,I,J,K)$ be a hypercomplex manifold and
\[
\Lambda^i_+(M) = \bigoplus_{p+q=i}\Lambda^{p,q}_{+,I}(M)
\]
the Hodge decomposition of qD-complex defined above.
Then there is a natural isomorphism
\begin{equation}\label{_qD_decompo_Equation_}
\Lambda^{p,q}_{+,I}(M)\cong \Lambda^{0,p+q}(M,I).
\end{equation}

{\bf Proof:} The following lemma is clear.

\hfill

%%%%%%%%%%%%%%%%%%%%%%%%%%%%%%%%%%%%%%%%%%%%%%%%
\lemma
Let $(M,I,J,K)$ be a hypercomplex manifold,
$\dim_{\Bbb H}M=n$, and $p$ an integer, $0\leq p\leq 2n$.
Then $\Lambda^{0,p}(M,I)\subset \Lambda^p(M)$
is pure of weight $p$.

\hfill

{\bf Proof:} Consider the operator 
$W_I:\;\Lambda^*(M) \arrow \Lambda^*(M)$
mapping a form $\eta\in \Lambda^{p,q}(M,I)$
to $\1(p-q)\eta$. Clearly, $W_I$ acts as a 
generator of $\goth u(1)$, with $\goth u(1)$
associated to $\rho_I:\; U(1)\arrow \End(\Lambda^*(M))$.
By \ref{_su(2)_action_explici_Lemma_},
$W_I\in \goth{su}(2)$, where the $\goth{su}(2)$-action 
on $\Lambda^*(M)$ is associated with the standard action
of $SU(2)$. Writing $\goth{su}(2)$ explicitly in terms
of generators $W_I$, $W_J$ and $W_K$, we find that $W_I$ 
generates a Cartan subalgebra of $\goth{su}(2)$
(indeed, the corresponding Lie group is
a maximal compact torus of $SU(2)$).
Since the Cartan algebra 
$\C \cdot W_I$ acts on $\Lambda^{p,0}(M,I)$ with 
weight $p$, the space $\Lambda^{p,0}(M,I)$
is of weight $\geq p$. On the other hand,
$\Lambda^p(M)$ is a representation of weight $\leq p$
(\ref{_weights_on_forms_Remark_}).
Therefore, $\Lambda^{p,0}(M,I)$ is pure of weight $p$.
\endproof

\hfill

%%%%%%%%%%%%%%%%%%%%%%%%%%%%%%%%%%%%%%%%%%
\remark\label{_Lambda^0,p_I,+_Remark_}
This argument also implies that $\Lambda^{0,p}(M,I)$
coinsides with
$\Lambda^{0,p}_{+,I}(M)\subset\Lambda^p_+(M)$
(here we consider $\Lambda^p_+(M)$ as a maximal
$SU(2)$-invariant subspace of weight $p$ in $\Lambda^p(M)$).

\hfill

Now, \ref{_qD_decompo_expli_Proposition_}
is implied by the general machinery of
$SU(2)$-representations. If $R$ is a finite-dimensional
$SU(2)$-representation of weight $\leq p$, the Cartan 
algebra action splits $R$ onto weight components 
$R= \oplus R_i$, $i=-p$, $-p+2$, ... $p-2$, $p$
the weights of the root $\1 W_I$
acting on $R_i$ as a multiplication 
by $i$. Moreover, if $R$ is pure of weight $p$,
then all spaces $R_i$ are naturally isomorphic,
with isomorphism provided by the $SU(2)$-action.

In the case $R=\Lambda^p_+(M)$, the decomposition $R= \oplus R_i$
is precisely the Hodge decomposition, hence the spaces
 $\Lambda^{p,q}_{+,I}(M)$ are naturally
isomorphic to for all $p, q\geq 0$ satisfying $p+q=i$.
We proved \ref{_qD_decompo_expli_Proposition_}.
\endproof

%%%%%%%%%%%%%%%%%%%%%%%%%%%%%%%%%%%%%%%%%%%%%%%%
\subsection{The Hodge decomposition on qD-complex:
an explicit construction}
%%%%%%%%%%%%%%%%%%%%%%%%%%%%%%%%%%%%%%%%%%%%%%%%

The isomorphism \eqref{_qD_decompo_Equation_}
can be made explicit, and also multiplicative, in the
following way. Let $\goth R$ be an irreducible
2-dimensional representation of $SU(2)$.
Clearly, any irreducible $SU(2)$-representation of
weight $p$ is isomorphic to $S^p\goth R$
(the $p$-th symmetric power of $\goth R$).
Consider the root  $\1 W_I\in \goth{su}(2)$,
constructed in Subsection \ref{_Hodge_on_qD_Subsection_}.
The corresponding $\goth{sl}(2)$-triple can be written as
\[
f= W_J+\1 W_K,\ \  \  g = W_J-\1 W_K, \ \ \  h= \1 W_I.
\]
Let $x, y$ be a basis in $\goth R$, such that
$hx=x$, $hy=-y$, $gx=y$, $fy=x$.

Consider a hypercomplex manifold $(M, I, J, K)$.
The bundle
\begin{equation}\label{_S^pR_otimes_Lambda^p_bundle_Equation_}
\goth S:= \bigoplus_p S^p{\goth R} \otimes \Lambda^{0,p}(M,I),
\end{equation}
is equipped with a natural multiplicative structure
(we assume that the elements of $S^p{\goth R}$
and $\Lambda^{0,q}(M,I)$ commute). We define
the following $SU(2)$-action on $\goth S$:
$SU(2)$ acts trivially on $\Lambda^{0,p}(M,I)$,
and in a standard way on $S^p{\goth R}$.

Consider an isomorphism 
${\goth R}\otimes \Lambda^{0,1}(M,I)\arrow\Lambda^1(M)$
mapping $x\otimes \eta$ to $J(\eta)$ and $y\otimes \eta$
to $\eta$. This map is clearly $SU(2)$-invariant.
Using the multiplicative structure on $\goth S$,
it can be extended to an $SU(2)$-invariant
algebra homomorphism
\begin{equation}\label{_S^pR_otimes_Lambda^p_to_+_Equation_}
 \bigoplus_p S^p{\goth R} \otimes \Lambda^{0,p}(M,I) \arrow\Lambda^*_+(M).
\end{equation}

\hfill

%%%%%%%%%%%%%%%%%%%%%%%%%%%%%%%%%%%%%%%%%%%%%%%%
\proposition \label{_qD_decompo_via_x,y_Proposition_}
In these assumptions, 
\eqref{_S^pR_otimes_Lambda^p_to_+_Equation_}
is an algebra isomorphism.

\hfill

{\bf Proof:} 
Let $\goth S^p\subset \goth S$ denote the grading $p$
component. Bijectivity of the map
\eqref{_S^pR_otimes_Lambda^p_to_+_Equation_}
is checked in the same way as one proves
\ref{_qD_decompo_expli_Proposition_}:
the Hodge components of $\goth S^p$ 
are all isomorphic, because $\goth S^p$
is a pure representation of weight $p$,
and the same is true for $\Lambda^p_+(M)$.
Therefore, it suffices to prove 
that the restriction of 
\eqref{_S^pR_otimes_Lambda^p_to_+_Equation_}
to one Hodge component, say, 
$y^p\Lambda^{0,p}(M,I)$, induces an
isomorphism
\[
y^p\Lambda^{0,p}(M,I)\arrow \Lambda^{0,p}_{+,I}(M).
\]
This is implied by the equality $\Lambda^{0,p}(M,I)=\Lambda^{0,p}_{+,I}(M)$
(\ref{_Lambda^0,p_I,+_Remark_}).
\endproof

%%%%%%%%%%%%%%%%%%%%%%%%%%%%%%%%%%%%%%%%%%%%%%
\subsection{The $\bar\6_J$-operator}
%%%%%%%%%%%%%%%%%%%%%%%%%%%%%%%%%%%%%%%%%%%%%%

Let $(M,I,J,K)$ be a hypercomplex manifold.
We extend \[ J:\; \Lambda^1(M) \arrow \Lambda^1(M)\]
to $\Lambda^*(M)$ by multiplicativity. Recall that 
\[ J(\Lambda^{p,q}(M,I))=\Lambda^{q,p}(M,I), \]
because $I$ and $J$ anticommute on $\Lambda^1(M)$.
Denote by 
\[ \bar\6_J:\;  \Lambda^{p,q}(M,I)\arrow \Lambda^{p,q+1}(M,I)
\]
the operator $J\circ \6 \circ J$, where
$\6:\;  \Lambda^{p,q}(M,I)\arrow \Lambda^{p+1,q}(M,I)$
is the standard Dolbeault operator on $(M,I)$, that is, the
$(1,0)$-part of the de Rham differential.
Since $\6^2=0$, we have $\bar\6_J^2=0$.
Since $I,J,K$ are integrable, the operators 
$d$, $d_I:= I\circ d\circ I$,  $d_J:= J\circ d\circ J$,  
$d_K:= K\circ d\circ K$ pairwise anticommute.
Therefore,
$\bar\6= \frac{d-\1d_I}{2}$ anticommutes with
$\bar\6_J= \frac{d_J-\1d_K}{2}$.
Writing the supercommutator as $\{\cdot, \cdot\}$,
we express this as
\begin{equation}\label{_commute_6_J_6_Equation_}
\{\bar\6_J, \bar\6_J \}=0, \ \ \ \{\bar\6_J, \bar\6 \}=0.
\end{equation}

%%%%%%%%%%%%%%%%%%%%%%%%%%%%%%%%%%%%%%%%%%%%%%
\subsection{The $\bar\6$, $\bar\6_J$-bicomplex}
\label{_6_6_J_bicomplex_Subsection_}
%%%%%%%%%%%%%%%%%%%%%%%%%%%%%%%%%%%%%%%%%%%%%%

Consider the quaternionic Dolbeault complex
$(\Lambda^*_+(M), d_+)$ constructed in Subsection
\ref{_qD_Subsection_}. Using the Hodge decomposition, we can
represent this complex as

\begin{equation}\label{_bicomple_d_+_Equation}
\begin{minipage}[m]{0.85\linewidth}
\begin{center}{ $
\xymatrix @C+1mm @R+10mm@!0  { 
  && \Lambda^0_{+,I}(M) \ar[dl]^{d^{1,0}_{+,I}} \ar[dr]^{d^{0,1}_{+,I}} 
   &&  \\
 & \Lambda^{1,0}_{+,I}(M) \ar[dl]^{d^{1,0}_{+,I}} \ar[dr]^{d^{0,1}_{+,I}} &
 & \Lambda^{0,1}_{+,I}(M) \ar[dl]^{d^{1,0}_{+,I}} \ar[dr]^{d^{0,1}_{+,I}}&\\
 \Lambda^{2,0}_{+,I}(M) && \Lambda^{1,1}_{+,I}(M) 
   && \Lambda^{0,2}_{+,I}(M) \\
}
$
}\end{center}
\end{minipage}
\end{equation}
where $d^{1,0}_{+,I}$,  $d^{0,1}_{+,I}$ are the Hodge components of 
the quaternionic Dolbeault differential $d_+$, taken with
respect to $I$. 

\hfill

Consider a hypercomplex manifold $(M, I, J, K)$.
Let 
\begin{equation}\label{_goth_S_iso_Equation_}
\bigoplus_p S^p{\goth R} \otimes \Lambda^{0,p}(M,I) \arrow\Lambda^*_+(M).
\end{equation}
be the isomorphism constructed in 
\ref{_qD_decompo_via_x,y_Proposition_}.
Writing the basis $x,y$ of $\goth R$ as in 
the proof of \ref{_qD_decompo_via_x,y_Proposition_},
we may write the Hodge decomposition of
\eqref{_goth_S_iso_Equation_} as
\[
x^py^q\Lambda^{0,p+q}(M,I)\cong \Lambda^{p,q}_{+,I}(M).
\]

%%%%%%%%%%%%%%%%%%%%%%%%%%%%%%%%%%%%%%%%%%%%%%
\theorem\label{_bico_ide_Theorem_}
Under this correspondence, $d^{0,1}_+$
corresponds to $\bar\6$ and $d^{1,0}_+$
to $\bar\6_J$. This way the 
bicomplex \eqref{_bicomple_d_+_Equation}
becomes equivalent to the 
bicomplex $(\Lambda^{0,p}(M,I),\bar\6, \bar\6_J)$
as follows:

\begin{equation}\label{_bicomple_XY_isomo_Equation}
\begin{minipage}[m]{0.85\linewidth}
{\tiny $
\xymatrix @C+1mm @R+10mm@!0  { 
  && \Lambda^0_{+,I}(M) \ar[dl]^{d^{1,0}_{+,I}} \ar[dr]^{d^{0,1}_{+,I}} 
   && && && \Lambda^{0,0}_I(M) \ar[dl]^{x \bar\6_J} \ar[dr]^{y \bar\6}
   &&  \\
 & \Lambda^{1,0}_{+,I}(M) \ar[dl]^{d^{1,0}_{+,I}} \ar[dr]^{d^{0,1}_{+,I}} &
 & \Lambda^{0,1}_{+,I}(M) \ar[dl]^{d^{1,0}_{+,I}} \ar[dr]^{d^{0,1}_{+,I}}&& 
\text{\large $\cong$} &
 &x\Lambda^{0,1}_I(M)\ar[dl]^{x \bar\6_J} \ar[dr]^{y \bar\6}&  &
 y\Lambda^{0,1}_I(M)\ar[dl]^{x \bar\6_J} \ar[dr]^{y \bar\6}&\\
 \Lambda^{0,2}_{+,I}(M) && \Lambda^{1,1}_{+,I}(M) 
   && \Lambda^{0,2}_{+,I}(M)& \ \ \ \ \ \ & x^2\Lambda^{0,2}_I(M)& & 
xy\Lambda^{0,2}_I(M) & & y^2\Lambda^{0,2}_I(M) \\
}
$
}
\end{minipage}
\end{equation}

{\bf Proof:} Consider the action of $x\bar\6_J+y\bar\6$
on 
\[ \bigoplus_p S^p{\goth R}\otimes \Lambda^{0,p}(M,I)\cong
   \Lambda^*_+(M)
\]
defined as in \eqref{_bicomple_XY_isomo_Equation}.
To prove \ref{_bico_ide_Theorem_}, we need to show that
\begin{equation}\label{_d_+_iso_6_J_Equation_}
x\bar\6_J+y\bar\6= d_+.
\end{equation}
Both of these operators satisfy the Leibniz rule,
hence it suffices to check \eqref{_d_+_iso_6_J_Equation_}
on some set of multiplicative generators of 
$\Lambda^*_+(M)$. On $\Lambda^0_+(M)$,
the equality \eqref{_d_+_iso_6_J_Equation_}
is clear from the definitions:
\begin{equation}\label{_d_+_iso_on_Lambda^0_Equation_}
x\bar\6_J+y\bar\6\restrict{\Lambda^0_+(M)}= d_+\restrict{\Lambda^0_+(M)}.
\end{equation}
 It is easy to check
that the space $\Lambda^0(M)\oplus d \Lambda^0(M)$
generates the algebra $\Lambda^*(M)$. Therefore,
$\Lambda^0_+(M)\oplus d_+ \Lambda^0_+(M)$
generates $\Lambda^*_+(M)$. To prove
\ref{_bico_ide_Theorem_}
it remains to show that
\begin{equation}\label{_d_+_iso_on_d_+_Lambda^0_Equation_}
x\bar\6_J+y\bar\6\restrict{d_+\Lambda^0_+(M)}= d_+\restrict{d_+\Lambda^0_+(M)}.
\end{equation}
Since $d_+^2=0$, $d_+\restrict{d_+\Lambda^0_+(M)}=0$.
By \eqref{_commute_6_J_6_Equation_},
\begin{equation}\label{_x_6_J+y_6_^2_Equation_}
(x\bar\6_J+y\bar\6)^2=0.
\end{equation}
 Using
\eqref{_d_+_iso_on_Lambda^0_Equation_} and 
\eqref{_x_6_J+y_6_^2_Equation_}, we obtain
\[
x\bar\6_J+y\bar\6\restrict{d_+\Lambda^0_+(M)}=
x\bar\6_J+y\bar\6\restrict{x\bar\6_J+y\bar\6(\Lambda^0_+(M))}=0.
\]
Therefore,
\[
x\bar\6_J+y\bar\6\restrict{d_+\Lambda^0_+(M)}= d_+\restrict{d_+\Lambda^0_+(M)}=0.
\]
This proves \eqref{_d_+_iso_on_d_+_Lambda^0_Equation_}.
\ref{_bico_ide_Theorem_} is proven.
\endproof

%%%%%%%%%%%%%%%%%%%%%%%%%%%%%%%%%%%%%%%%%%%%%%

\section{Kodaira identities for qD-complex}

%%%%%%%%%%%%%%%%%%%%%%%%%%%%%%%%%%%%%%%%%%%%%%

%%%%%%%%%%%%%%%%%%%%%%%%%%%%%%%%%%%%%%%%%%%%%%%%%%%%%%%%%%%%%%
\subsection{The Lefschetz-type 
$\goth{sl}(2)$-action on
  $\Lambda^{0,*}(M,I)\otimes\End(B)$}
\label{_Lefschetz_Subsection_}
%%%%%%%%%%%%%%%%%%%%%%%%%%%%%%%%%%%%%%%%%%%%%%%%%%%%%%%%%%%%%%

Let $(M,I,J,K)$ be a hyperk\"ahler manifold, 
$B$ a holomorphic Hermitian vector bundle on $(M,I)$,
and  $\Lambda^{0,*}(M,I)\otimes B$
the space of $(0,p)$-forms with coefficients in
$B$. Denote by $\bar\Omega\in \Lambda^{0,2}(M,I)$
the standard $(0,2)$-form $\omega_J+\1\omega_K$
(Subsection \ref{_HK_intro_Subsection_}).

Using a hyperk\"ahler metric, we constract
a natural Hermitian structure on $\Lambda^{0,*}(M,I)\otimes B$. Denote 
by $L_{\bar \Omega}:\; \Lambda^{q,p}(M,I) \arrow
\Lambda^{q,p+2}(M)$ the operator of exterior
multiplication by $\bar\Omega$, and by
$\Lambda_{\bar \Omega}:\; \Lambda^{q,p}(M,I)\otimes B \arrow
\Lambda^{q,p-2}(M)\otimes B$ its Hermitian adjoint. The same
argument as proves the usual Lefschetz Theorem 
about the $\goth{sl}(2)$-action (see
\cite{_Griffi_Harri_}) can be used to prove
the following linear-algebraic result, which is due to A. Fujiki.

\hfill

%%%%%%%%%%%%%%%%%%%%%%%%%%%%%%%%%%%%%%%%%%%%%%%%%%%%%%%%%%%%%%
\proposition\label{_SL(2)_Proposition_}
(\cite{_Fujiki_}) In the above assumptions, let 
\[ 
  H_{\bar\Omega}:= [L_{\bar \Omega}, \Lambda_{\bar
  \Omega}]
\]
be a commutator of $L_{\bar \Omega}, \Lambda_{\bar
  \Omega}$. Then $H_{\bar\Omega}$ is a scalar operator,
multiplying a $(q,p)$-form by $n-p$, where $n=\frac 1 2
  \dim{\Bbb H}(M)$. Moreover,  $L_{\bar \Omega}, \Lambda_{\bar
  \Omega}$, $H_{\bar\Omega}$ is an $\goth{sl}(2)$-triple.

\hfill

{\bf Proof:} 
See \cite{_Verbitsky:Hyperholo_bundles_},
Theorem 4.2). \endproof

\hfill

Let $\theta\in \Lambda^{0,1}(M,I)\otimes\End(B)$ be
a 1-form. Denote by 
\[ L_\theta:\; \Lambda^{q,p}(M,I)\otimes B \arrow
\Lambda^{q,p+1}(M,I)\otimes B \] the operator of
multiplication by $\theta$, and 
let 
\[ \Lambda_\theta:\; 
   \Lambda^{q,p}(M,I) \arrow
   \Lambda^{q,p-1}(M,I)
\]
be its Hermitian adjoint. Denote by $\theta_J$ the
$(0,1)$-form $J(\bar\theta)$.

\hfill

%%%%%%%%%%%%%%%%%%%%%%%%%%%%%%%%%%%%%%%%%%%%%%%%%%%%%%%
\claim\label{_L_Omega_with_Lambda_theta_Claim_}
In the above assumptions, we have
\begin{equation}\label{_[L_Omega_Lambda_theta_]_Equation_}
[ L_{\bar\Omega}, \Lambda_\theta] = L_{\theta_J}.
\end{equation}

{\bf Proof:} Follows
from a trivial computation.
\endproof

%%%%%%%%%%%%%%%%%%%%%%%%%%%%%%%%%%%%%%%%%%%%%%%%%%%%%%%%%%%%%%
\subsection{$\bar\6$, $\bar\6_J$ with coefficients in a
  bundle}
\label{_6,6_J_in_bundle_Subsection_}
%%%%%%%%%%%%%%%%%%%%%%%%%%%%%%%%%%%%%%%%%%%%%%%%%%%%%%%%%%%%%%

Let $(M,I,J,K)$ be a hyperk\"ahler manifold, and
$B$ a holomorphic Hermitian vector bundle on $(M,I)$.
Consider the standard (Chern) Hermitian connection 
$\nabla$ on $B$, $\nabla= \nabla^{1,0}+ \bar\6$, where
$\bar\6:\; B \arrow B \otimes \Lambda^{0,1}(M,I)$
is the holomorphic structure operator. Denote by 
$\bar\6_J:\; B \arrow B \otimes \Lambda^{0,1}(M,I)$
the composition of $\nabla^{1,0}:\; B \otimes
\Lambda^{1,0}(M,I)$ and  
\[ \Id_B\otimes J:\;  B \otimes \Lambda^{1,0}(M,I)\arrow B \otimes
   \Lambda^{0,1}(M,I)
\]
be an endomorphism associated with $J\in {\Bbb H}$.
We extend $\bar\6$, $\bar\6_J$
to operators
\[
\bar\6, \bar\6_J:\; \Lambda^{0,p}(M,I)\otimes B \arrow
   \Lambda^{0,p+1}(M,I)\otimes B.
\]
using the Leibniz rule.

\hfill

%%%%%%%%%%%%%%%%%%%%%%%%%%%%%%%%%%%%%%%%%%%%%%%%%%%%%%%%%%%%%%
\proposition\label{_6_6_J_commu_on_bundles_Proposition_}
In these assumptions. $\bar\6^2=\bar\6_J^2=0$, and
the anticommutator $\{\bar\6, \bar\6_J\}$
acts on $\Lambda^{0,*}(M,I)$ as a multiplication
by an $\End(B)$-valued 2-form 
$\Theta_+\in \Lambda^{0,2}(M,I)\otimes \End B$.
Moreover, under the identification
\[
\Lambda^{0,2}(M,I)\otimes \End B\cong \Lambda^{1,1}_{+,I}(M)\otimes \End B
\]
(\ref{_qD_decompo_expli_Proposition_}),
$\Theta_+$ corresponds to the 
$\Lambda^2_+(M)$-part
of the curvature of $B$.

\hfill

{\bf Proof:} Let 
\[ \nabla_+:\; \Lambda^p_+(M)\otimes
B\arrow\Lambda^{p+1}_+(M)\otimes B
\]
be the connection operator restricted to
$\Lambda^*_+(M)\otimes B$, and
$\nabla_+=\nabla^{1,0}_+ + \6_+$ its Hodge decomposition.
Clearly, $\nabla_+^2$ is the $\Lambda^2_+(M)$-part
of the curvature of $B$.

Now, \ref{_6_6_J_commu_on_bundles_Proposition_}
follows immediately from 
\ref{_bico_ide_Theorem_}. Indeed, under the isomorphism
\eqref{_bicomple_XY_isomo_Equation}, $x\bar\6_J$
corresponds to $\nabla^{1,0}_+$; since the curvature
of the Chern connection is of type $(1,1)$, we have
$(\nabla^{1,0}_+)^2=0$, hence $\bar\6_J^2=0$. 
Similarly, the operator $\{x\bar\6_J, y \bar\6\}$ under 
the isomorphism \eqref{_bicomple_XY_isomo_Equation}
corresponds to $\{ \nabla^{1,0}_+, \6_+\} = \nabla_+^2$.
\endproof

%%%%%%%%%%%%%%%%%%%%%%%%%%%%%%%%%%%%%%%%%%%%%%%%%%%%%%%%%%%%%%
\subsection{Kodaira relations for $\bar\6$, $\bar\6_J$}
%%%%%%%%%%%%%%%%%%%%%%%%%%%%%%%%%%%%%%%%%%%%%%%%%%%%%%%%%%%%%%

Let $(M, I, J, K)$ be a hyperk\"ahler manifold.
Consider the bicomplex 
\[ (\Lambda^{0,*}(M,I), \bar\6,
\bar\6_J), 
\]
constructed in Subsection
\ref{_6_6_J_bicomplex_Subsection_}.
Let 
\[  L_{\bar \Omega}:\;
\Lambda^{0,p}(M,I)\arrow\Lambda^{0,p+2}(M,I)
\]
be an operator of exterrior multiplication by
$\bar\Omega$ (Subsection \ref{_Lefschetz_Subsection_}),
and 
\[ \bar\6^*, \bar\6^*_J:\;
  \Lambda^{0,p}(M,I)\arrow\Lambda^{0,p-1}(M,I).
\]
the operators Hermitian
adjoint to $\bar\6$, $\bar\6_J$. 

\hfill

The following proposition is well known.

\hfill

%%%%%%%%%%%%%%%%%%%%%%%%%%%%%%%%%%%%%%%%%%%%%%%%%%%%%%%%%%%%%%
\proposition\label{_Kodaira_Proposition_}
In these assimptions, the following commutator relations
hold.
\begin{equation}\label{_Kodaira_Equation_}
[L_{\bar \Omega}, \bar\6^*] =  \bar\6_J, \ \ 
[L_{\bar \Omega}, \bar\6_J^*] = - \bar\6.
\end{equation}

{\bf Proof:} The proof of \eqref{_Kodaira_Equation_}
is essentially the same as the proof of the usual
Kodaira relations; see
e.g. \cite{_Verbitsky:Hyperholo_bundles_}, Proposition
4.2. \endproof

\hfill

The same argument, applied locally to $\End(B)$-valued
forms, gives the following theorem.

\hfill

%%%%%%%%%%%%%%%%%%%%%%%%%%%%%%%%%%%%%%%%%%%%%%%%%%%%%%
\theorem\label{_Kodaira_main_Theorem_}
Let $(M, I, J, K)$ be a hyperk\"ahler manifold,  $B$ 
a holomorphic Hermitian vector bundle on $(M,I)$, 
\[
\bar\6, \bar\6_J:\; \Lambda^{0,p}(M,I)\otimes B \arrow
   \Lambda^{0,p+1}(M,I)\otimes B.
\]
the operators constructed in Subsection
\ref{_6,6_J_in_bundle_Subsection_},
and $ \bar\6^*, \bar\6_J^*$ the Hermitian
adjoint operators. Then 
\begin{equation}\label{_Kodaira_in_bundle_Equation_}
[L_{\bar \Omega}, \bar\6^*] = \bar\6_J, \ \ 
[L_{\bar \Omega}, \bar\6_J^*] = - \bar\6.
\end{equation}

{\bf Proof:} See \cite{_Verbitsky:Hyperholo_bundles_}.
\endproof
\endproof

%%%%%%%%%%%%%%%%%%%%%%%%%%%%%%%%%%%%%%%%%%%%%%%%%%%%%%%%%%%%%%
\subsection{Kodaira-Nakano identities}
%%%%%%%%%%%%%%%%%%%%%%%%%%%%%%%%%%%%%%%%%%%%%%%%%%%%%%%%%%%%%%

The following theorem is the qD-analogue of the usual
Kodaira-Nakano identity (or, rather, the identity
used in the proof of Kodaira-Nakano vanishing)

\hfill

%%%%%%%%%%%%%%%%%%%%%%%%%%%%%%%%%%%%%%%%%%%%%%%%%%%%%%
\theorem\label{_Kodaira_Nakano_Laplacians_Theorem_}
Let $(M, I, J, K)$ be a hyperk\"ahler manifold,  $B$ 
a holomorphic Hermitian vector bundle on $(M,I)$,
\[
\bar\6, \bar\6_J:\; \Lambda^{0,p}(M,I)\otimes B \arrow
   \Lambda^{0,p+1}(M,I)\otimes B.
\]
the operators constructed in Subsection
\ref{_6,6_J_in_bundle_Subsection_},
and $ \bar\6^*, \bar\6_J^*$ the Hermitian
adjoint operators. Consider the Laplacians
\[
\Delta_{\bar\6}:= \{ \bar\6, \bar\6^*\}, \ \ 
\Delta_{\bar\6_J}:= \{ \bar\6_J, \bar\6_J^*\}
\]
(here, as elsewhere, $\{\cdot, \cdot\}$ denotes the
anticommutator).
Then
\begin{equation}\label{_Laplacians_K-N_ide_Equation_}
\Delta_{\bar\6}-\Delta_{\bar\6_J}=  [ \Theta_+, \Lambda_{\bar\Omega}],
\end{equation}
where 
\[ \Theta_+:\; \Lambda^{0,p}(M,I)\otimes B \arrow
   \Lambda^{0,p+2}(M,I)
\]
is an operator defined as 
\[
\Theta_+:=\{\bar\6, \bar\6_J\}
\]
and identified with the
$\Lambda_+^2(M)\otimes \End B$-part of the curvature  of $B$
as in \ref{_6_6_J_commu_on_bundles_Proposition_}.

\hfill

{\bf Proof:} Using the graded Jacobi identity and 
\ref{_Kodaira_main_Theorem_}, we obtain
\[
[ \Theta_+, \Lambda_{\bar\Omega}]=
- [ \Lambda_{\bar\Omega}, \{\bar\6, \bar\6_J\}] =
\{ \bar\6, \bar\6^*\} - \{ \bar\6_J, \bar\6_J^*\}= 
\Delta_{\bar\6}-\Delta_{\bar\6_J}.
\]
\endproof

%%%%%%%%%%%%%%%%%%%%%%%%%%%%%%%%%%%%%%%%%%%%%%%%%%%%%%%%%%%%%%

\section{Cohomology of hyperk\"ahler manifolds}
\label{_BBF_Section_}

%%%%%%%%%%%%%%%%%%%%%%%%%%%%%%%%%%%%%%%%%%%%%%%%%%%%%%%%%%%%%%

For the convenience of the reader, we recall here
some well-known facts about the structure of $H^2(M)$
for $M$ a compact, irreducible hyperk\"ahler manifold; see
\cite{_Bogomolov_}, \cite{_Besse:Einst_Manifo_}, 
\cite{_Beauville_} and \cite{_Fujiki_} for details.

%%%%%%%%%%%%%%%%%%%%%%%%%%%%%%%%%%%%%%%%%%%%%%%%%%%%%%%%%%%%%%
\subsection{$SU(2)$-action on $H^2(M)$}
\label{_SU(2)_on_coho_Subsection_}
%%%%%%%%%%%%%%%%%%%%%%%%%%%%%%%%%%%%%%%%%%%%%%%%%%%%%%%%%%%%%%

Let $(M, I, J, K, g)$ be a compact, irreducible hyperk\"ahler
manifold. Since $g$ is K\"ahler with respect to
$(I,J,K)$, we have
\[
\nabla I = \nabla J = \nabla K=0,
\]
where $\nabla$ denotes the Levi-Civita connection.
Chern has shown that covariantly constant endomorphisms
of $\Lambda^*(M)$ commute with the Laplacian 
(see \cite{_Besse:Einst_Manifo_}). Then the
$SU(2)$-action generated by $I, J, K\in {\Bbb H}^*$
also commutes with the Laplacian. This gives an 
$SU(2)$-action on the space  of harmonic forms on $M$.
Identifying the harmonic forms with cohomology,
we obtain an $SU(2)$-action on the cohomology as well.

Let $H^2(M) = H^2_+(M)\oplus H^2(M)_{SU(2)-inv}$
be a decomposition of $H^2(M)$ onto its weight 2
and weight 0 components.
Using the weights of the Cartan algebra action
as in the proof of \ref{_qD_decompo_expli_Proposition_},
we find that 
\[
\dim H^{2,0}(M,I) = \dim H^{1,1}_+(M,I)=\dim H^{0,2}(M,I).
\]
Since $M$ is irreducible, $\dim H^{2,0}(M,I)=1$
and the space $H^{1,1}_+(M,I)$ is one-dimensional.
Let $H^2(M)_{SU(2)-inv}$ be 
the space of $SU(2)$-invariant classes.
It is easy to check that $SU(2)$-invariant classes are 
all of type $(1,1)$ (e.g. \cite{_Verbitsky:Hyperholo_bundles_}).

Since $H^{1,1}_+(M,I)$ is one-dimensional and
generated by the K\"ahler form $\omega_I$, we have
a decomposition
\begin{equation}\label{_SU_2_on_H^1,1_Equation_}
H^{1,1}(M,I)=\C \omega_I\oplus H^2(M)_{SU(2)-inv}.
\end{equation}

Using the $\goth{so}(1,4)$-action generated
by the three Lefschetz $\goth{sl}(2)$-triples
associated with the K\"ahler structures 
$I, J, K$ as in \cite{_so(5)_},
we can easily show that an $SU(2)$-invariant
2-form is primitive\footnote{Recall that the {\bf primitive 
classes} (\cite{_Griffi_Harri_})
are cohomology classes which
  satisfy $\Lambda(\eta)=0$, where $\Lambda:\;
  H^i(M)\arrow H^{i-2}(M)$ is the dual
Lefschetz operator. A $(1,1)$-class is primitive if and only if it is
orthogonal to the K\"ahler form with respect to the
Riemann-Hodge pairing.} (see
e.g. \cite{_Verbitsky:Hyperholo_bundles_}).

\hfill

This gives the following well-known statement
(\cite{_Verbitsky:Hyperholo_bundles_}).

\hfill

%%%%%%%%%%%%%%%%%%%%%%%%%%%%%%%%%%%%%%%%%%%%%%%%%%
\claim\label{_primi_SU(2)_Claim_}
Let $(M, I, J, K)$ be a compact, irreducible hyperk\"ahler
manifold. Then the space  $H^{1,1}_{prim}(M,I)$ of primitive
classes in $H^{1,1}(M,I)$ coincides with the
space $H^2(M)_{SU(2)-inv}$ of $SU(2)$-invariant classes.

\hfill

{\bf Proof:} Since all $SU(2)$-invariant
classes are primitive, 
$H^{1,1}_{prim}(M,I)$ contains $H^2(M)_{SU(2)-inv}$.
Comparing the decomposition
\eqref{_SU_2_on_H^1,1_Equation_} with
\[
H^{1,1}(M,I) = H^{1,1}_{prim}(M,I)\oplus \C\omega_I,
\]
we find that $\dim H^{1,1}_{prim}(M,I)=\dim H^2(M)_{SU(2)-inv}$.
\endproof

%%%%%%%%%%%%%%%%%%%%%%%%%%%%%%%%%%%%%%%%%%%%%%%%%%%%%%%%%%%%%%
\subsection{The Bogomolov-Beauville-Fujiki form}
%%%%%%%%%%%%%%%%%%%%%%%%%%%%%%%%%%%%%%%%%%%%%%%%%%%%%%%%%%%%%%

Let $(M,I,J,K)$ be a compact hyperk\"ahler manifold,
and $\Omega:= \omega_J+\1\omega_K$ the holomorphic
symplectic form on $(M,I)$. F. Bogomolov 
(\cite{_Bogomolov_}) defined the
following bilinear symmetric 2-form on $H^{1,1}(M,I)$:
\begin{equation}\label{_BBF_form_on_H^11_Equation_}
\tilde q(\eta,\eta'):= \int_M\eta\wedge \eta'\wedge \Omega^{n-1}
\wedge \bar\Omega^{n-1},
\end{equation}
where $n=\dim{\Bbb H}M$.
Since $\Omega\wedge\bar\Omega$ is a positive (2,2)-form,
$\tilde q$ is positive on the K\"ahler cone of $(M,I)$: 
\begin{equation}\label{_q_on_K_cone_Equation_}
\forall \omega\in {\cal K}\ \ \ \tilde q(\omega, \omega)>0
\end{equation}
An elementary
linear-algebraic calculation similar to the proof 
of Riemann-Hodge bilinear relations implies that
$\tilde q(\eta,\eta)<0$ for $\eta$ 
primitive.
Therefore, $\tilde q$ has signature $(+, -, -, -, \dots)$
on $H^{1,1}(M,I)\cap H^2(M, \R)$.

The form $\tilde q$ is topological by its nature.

\hfill

%%%%%%%%%%%%%%%%%%%%%%%%%%%%%%%%%%%%%%%%%%%%%%%%
\theorem\label{_BBF_Theorem_}
(\cite{_Fujiki_})
Let $(M, I, J, K)$ be a compact, irreducible hyperk\"ahler
manifold of real dimension $4n$. Then there exists a
bilinear, symmetric non-degenerate
2-form $q:\; H^2(M,\Q)\otimes H^2(M,\Q)\arrow \Q$
such that 
\begin{equation}\label{_BBF_form_Equation_}
\int_M \eta^{2n}= q(\eta,\eta)^n,
\end{equation}
for all $\eta\in H^2(M)$. Moreover,
$q$ is proportional to the form 
\eqref{_BBF_form_on_H^11_Equation_} on $H^{1,1}(M)$, and has signature
$(+,+,+, -,-,-, ...)$.

\endproof

\hfill

%%%%%%%%%%%%%%%%%%%%%%%%%%%%%%%
\remark
If $n$ is odd, the equation \eqref{_BBF_form_Equation_}
determines $q$ uniquely, otherwise -- up to a sign.
To choose a sign, we use \eqref{_q_on_K_cone_Equation_}

\hfill

%%%%%%%%%%%%%%%%%%%%%%%%%%%%%%%%%%%%%%%%%%%%%%%
\definition\label{_BBF_Definition_}
Let $(M, I, J, K)$ be a compact, irreducible hyperk\"ahler
manifold. A {\bf Beauville-Bogomolov-Fujiki form} on $M$
is a form $q:\; H^2(M,\Q)\otimes H^2(M,\Q)\arrow \Q$
which satisfies \eqref{_BBF_form_Equation_},
and take positive values on the K\"ahler cone of $(M,I)$.
Such a form always exists and is unique, by \ref{_BBF_Theorem_}.

\hfill

%%%%%%%%%%%%%%%%%%%%%%%%%%%%%%%
\remark
The  Beauville-Bogomolov-Fujiki form is integer,
but not unimodular on $H^2(M, \Z)$.

\hfill
 
The Beauville-Bogomolov-Fujiki form can be expressed 
in terms of the $SU(2)$-\-ac\-tion on cohomology 
(Subsection \ref{_SU(2)_on_coho_Subsection_}) as follows:

\hfill

%%%%%%%%%%%%%%%%%%%%%%%%%%%%%%%%%%%%%%%%
\claim\label{_BBF_via_SU(2)_Claim_}
Let $(M,I,J,K)$ be a compact, irreducible hyperk\"ahler
manifold, and $(\cdot,\cdot)_{\cal H}$ the positive
definite pairing on cohomology associated with the
Euclidean metric on the space of harmonic forms
induced by the Riemannian structure.
Consider the form $q'$ which is equal to 
$(\cdot,\cdot)_{\cal H}$ on the 3-dimensional
space generated by $\omega_I$, $\omega_J$, $\omega_K$,
and to $-(\cdot,\cdot)_{\cal H}$ on its orthogonal
complement. Then $q'$ is proportional to the 
Beauville-Bogomolov-Fujiki form.

\hfill

{\bf Proof:} See e.g. \cite{_coho_announce_},
Theorem 2.1. \endproof

\hfill

This immediately gives the following corollary.

\hfill

%%%%%%%%%%%%%%%%%%%%%%%%%%%%%%%%%%%%%%%%%%%%%%%%%
\corollary\label{_q_SU(2)_inv_Corollary_}
Consider the natural $SU(2)$-action on the cohomology 
of a hyperk\"ahler manifold. Then the
Beauville-Bogomolov-Fujiki form is $SU(2)$-invariant.

\endproof

\hfill

Using the Hodge-Riemann bilinear relations, we can express
$(\cdot,\cdot)_{\cal H}$ in terms of the product structure on
cohomology. Together with \ref{_BBF_via_SU(2)_Claim_},
this gives 
\begin{equation}\label{_q_via_HR_Equation_}
q'(\eta_1,\eta_2) =   \int_X \omega_I^{2n-2}\wedge \eta_1\wedge\eta_2  -
   \frac{2n-2}{(2n-1)^2}\cdot 
   \frac{\int_X \omega_I^{2n-1}\eta_1 \cdot \int_X\omega_I^{2n-1}\eta_2}
   {\int_X \omega_I^n}
\end{equation}
for any $\eta_1, \eta_2\in H^2(M)$
(see \cite{_Verbitsky:cohomo_}, Claim 5.1).
This formula is due to 
A. Beauville.

\hfill

The following claim follows directly from 
\eqref{_q_via_HR_Equation_} and \ref{_primi_SU(2)_Claim_}.

\hfill

%%%%%%%%%%%%%%%%%%%%%%%%%%%%%%%%%%%%%%%%%%%%%%%%%%
\claim\label{_BBF_and_primitive_Claim_}
Let $(M,I,J,K)$ be a irreducible, compact hyperk\"ahler manifold,
and $\eta \in H^{1,1}(M,I)$ a $(1,1)$-class. Then the 
following assertions are equivalent.
\begin{description}
\item[(i)] $q(\eta,\omega_I)=0$, where $q$ is Beauville-Bogomolov-Fujiki
form, and $\omega_I$ the K\"ahler class of $(M,I)$
\item[(ii)] $\eta$ is primitive
\item[(iii)] $\eta$ is $SU(2)$-invariant.
\end{description}
{\bf Proof:} The equivalence of (ii) and (iii) is implied
by \ref{_primi_SU(2)_Claim_}, and the equivalence of
(i) and (iii) by \eqref{_q_via_HR_Equation_}.
\endproof

%%%%%%%%%%%%%%%%%%%%%%%%%%%%%%%%%%%%%%%%%%%%%%%%%%%%%%%%%%%%%%

\section{The vanishing of cohomology}

%%%%%%%%%%%%%%%%%%%%%%%%%%%%%%%%%%%%%%%%%%%%%%%%%%%%%%%%%%%%%%

%%%%%%%%%%%%%%%%%%%%%%%%%%%%%%%%%%%%%%%%%%%%%%%%%%%%%%%%%%%%%%
\subsection[Cohomology vanishing for  line bundles with
  $q(c_1(L), \omega)>0$]{Cohomology vanishing for  line bundles with\\
  $q(c_1(L), \omega)>0$}
%%%%%%%%%%%%%%%%%%%%%%%%%%%%%%%%%%%%%%%%%%%%%%%%%%%%%%%%%%%%%%

The following result is implied immediately by the
quaternionic Kodaira-Nakano identity
(\ref{_Kodaira_Nakano_Laplacians_Theorem_}),
in the same fashion as the usual Kodaira-Nakano vanishing
follows from the the usual Kodaira-Nakano identity.

\hfill

%%%%%%%%%%%%%%%%%%%%%%%%%%%%%%%%%%%%%%%%%%%%%%%%%%%%%%%%%%%%%%
\proposition\label{_vanishing_for_q(L,omega)>0_Proposition_}
Let $(M,I,J,K)$ be a compact, irreducible hyperk\"ahler
manifold, $\dim_{\Bbb R}M=4n$, and
$L$ a holomorphic line bundle on $(M,I)$,
such that $q(c_1(L), \omega_I)>0$, where $\omega_I$
is the K\"ahler class of $(M,I)$. Then the 
holomorphic cohomology $H^i(M,L)$ are zero for $i>n$.

\hfill

{\bf Proof:} Let $\eta$ be a harmonic form representing
$c_1(L)$. We may chose the Hermitian structure
on $L$ in such a way that $\eta$ is equal to the
curvature of $L$ (see \cite{_Griffi_Harri_}).
Let $\omega$ denote the K\"ahler form
of $(M,I)$. Abusing the notation,
we denote the K\"ahler class of $(M,I)$
by the same letter.

The cohomology class
\[
\kappa:= c_1(L) - \frac{q(c_1(L),\omega)}{q(\omega,\omega)}\omega
\]
clearly satisfies $q(\kappa,\omega)=0$. Therefore, $\kappa$ is
$SU(2)$-invariant (\ref{_BBF_and_primitive_Claim_}). 
Since $\eta$ is harmonic, the harmonic form
\begin{equation}\label{_primiti_decompo_H^2_Equation_}
\eta - \frac{q(c_1(L),\omega)}{q(\omega,\omega)}\omega
\end{equation}
representing $\kappa$ is also $SU(2)$-invariant. 
Let $\tilde \omega$ be the form
\[ \omega\in\Lambda^{1,1}_{+,I}(M)\]  considered as an
element in $\Lambda^{0,2}(M,I)$ using the isomorphism,
constructed \ref{_qD_decompo_expli_Proposition_}.
By \ref{_6_6_J_commu_on_bundles_Proposition_},
\begin{equation}\label{_Theta_+_for_line_bu_kahle_Equation_}
\Theta_+ = \{\bar\6,\bar\6_J\} = \frac{q(c_1(L),\omega)}{q(\omega,\omega)}\tilde\omega.
\end{equation}
Clearly, $\omega_I, \omega_J, \omega_K$ form
a 3-dimensional irreducible $SU(2)$-invariant 
subspace of $\Lambda^2(M)$. A trivial
calculation is used to show that
$\tilde\omega$ is in fact equal to $\bar\Omega$.
This gives
\begin{equation}\label{_Theta_+_via_Omega_Equation_}
\Theta_+ = \{\bar\6,\bar\6_J\} =\lambda L_{\bar\Omega},
\end{equation}
where $\lambda=\frac{q(c_1(L),\omega)}{q(\omega,\omega)}$
is a positive constant. Comparing
\eqref{_Theta_+_via_Omega_Equation_}, 
Kodaira-Nakano identity \eqref{_Laplacians_K-N_ide_Equation_}
and the quaternionic Lefschetz theorem (\ref{_SL(2)_Proposition_}),
we obtain
\begin{equation}\label{_H_as_commu_and_Lapla_Equation_}
\Delta_{\bar\6}-\Delta_{\bar\6_J}= [ \Theta_+,
  \Lambda_{\bar\Omega}] = \lambda H_{\bar\Omega}.
\end{equation}
On $(0,i)$-forms this operator acts as $(i-n)\lambda$.
Given a harmonic form
\[ \nu\in \ker \Delta_{\bar\6}\subset
  \Lambda^{0,i}(M,I)\otimes L,
\]
we can obtain
\begin{equation}\label{_Laplacian_on_harmo_via_H_Equation_}
0 = \Delta_{\bar\6}(\nu) = \Delta_{\bar\6_J}(\nu)+ (i-n)\lambda\nu.
\end{equation}
Since 
\begin{equation}\label{_positi_used_for_harmo_Equation_}
   (\Delta_{\bar\6_J}(\nu), \nu) = (\bar\6_J\eta,
   \bar\6_J\eta) +(\bar\6_J^*\eta,
   \bar\6_J^*\eta) \geq 0,
\end{equation}
\eqref{_Laplacian_on_harmo_via_H_Equation_}
leads to $(\nu, \nu)=0$ for $i>n$.
The harmonic $(0,i)$-forms are identified with the
$i$-th holomorphic cohomology of $L$ as usual.
We proved \ref{_vanishing_for_q(L,omega)>0_Proposition_}. \endproof

\hfill

%%%%%%%%%%%%%%%%%%%%%%%%%%%%%%%%%%%%%%%%%%%%%%%%%%%%%%%%%%%%%%
\remark \label{_positivity_Lapla_Remark_}
Let $W$ be a Hermitian vector space.
A {\bf positive operator} $A:\; W \arrow W$  is an operator which 
satisfies $(A(x),x)\geq 0$ for all $x\in W$. $A$ is 
{\bf positive definite} if this inequality is
strict for all non-zero $x$. From 
\eqref{_positi_used_for_harmo_Equation_},
we obtain that the Laplacians $\Delta_{\bar\6_J}$ and
$\Delta_{\bar\6}$ are positive. If 
\[ \Delta_{\bar\6} = \Delta_{\bar\6_J}+A,
\]
where $A$ is positive definite, then $\ker \Delta_{\bar\6}=0$.
This argument is used quite often in geometry and analysis.

\hfill

%%%%%%%%%%%%%%%%%%%%%%%%%%%%%%%%%%%%%%%%%%%%%%%%%%%%%%%%%%%%%%
\remark \label{_vanishing_for_q(L,omega)<0_Serre_Remark_}
Let $(M,I,J,K)$ be a compact, irreducible hyperk\"ahler
manifold, $\dim_{\Bbb H}M=n$, and 
$L$ a holomorphic line bundle on $(M,I)$.
Then Serre's duality gives $H^i(M,L)^*\cong H^{n-i}(M,L^*)$,
because the canonical class of $M$ is trivial.
Therefore, \ref{_vanishing_for_q(L,omega)>0_Proposition_}
implies that $H^i(M,L)$ vanish for all $i<n$
if $L$ is a holomorphic line bundle
on $(M,I)$  with $q(c_1(L), \omega)<0$.

%%%%%%%%%%%%%%%%%%%%%%%%%%%%%%%%%%%%%%%%%%%%%%%%%%%%%%%%%%%%%%
\subsection{The dual K\"ahler cone and vanishing}
\label{_dual_Kahler_Subsection_}
%%%%%%%%%%%%%%%%%%%%%%%%%%%%%%%%%%%%%%%%%%%%%%%%%%%%%%%%%%%%%%

Let $(M,I)$ be a K\"ahler manifold. 

\hfill

%%%%%%%%%%%%%%%%%%%%%%%%%%%%%%%%%%%%%%%%%%%%%%%
\definition \label{_Kahler_cone_Definition_}
The {\bf K\"ahler cone} ${\cal K}\subset H^{1,1}(M,I)$
for $(M,I)$ is the set of all K\"ahler classes
$\omega\in H^{1,1}_I(M, \R)$, where $H^{1,1}_I(M, \R)$
denotes the intersection $H^{1,1}(M,I)\cap H^2(M, \R)$.
Clearly, ${\cal K}$ is a convex cone in $H^{1,1}_I(M, \R)$.

\hfill

Now, let $(M,I,J,K)$ be a hyperk\"ahler manifold,
and ${\cal K}\subset H^{1,1}_I(M, \R)$ the K\"ahler
cone of $(M,I)$, and 
$q:\; H^{1,1}_I(M, \R)\times H^{1,1}_I(M, \R)\arrow \R$
the Beauville-Bogomolov-Fujiki form.  We define
{\bf the dual K\"ahler cone} \[ {\cal K}\check{\;}\subset H^{1,1}_I(M, \R)\] as 
\[
{\cal K}\check{\;}:= \{ x\in H^{1,1}_I(M, \R) \ \ | \ \ \forall
y\in {\cal K},\ \  q(x,y)> 0 \}
\]
It is an open, convex cone. Since a product of
two K\"ahler forms is positive, we have 
${\cal  K}\check{\;}\supset {\cal K}$.

\hfill

Denote by $\bar{\cal  K}\check{\;}$ the closure of
${\cal  K}\check{\;}$  in $H^{1,1}_I(M, \R)$, and by
$-\bar{\cal  K}\check{\;}$ the opposite cone. Clearly,
\[ 
\bar{\cal K}\check{\;}:= \{ x\in H^{1,1}_I(M, \R) \ \ | \ \ \forall
y\in {\cal K},\ \  q(x,y)\geq 0 \}
\]
and
\[ 
-\bar{\cal K}\check{\;}:= \{ x\in H^{1,1}_I(M, \R) \ \ | \ \ \forall
y\in {\cal K},\ \  q(x,y)\leq 0 \}
\]

\hfill

\ref{_vanishing_for_q(L,omega)>0_Proposition_} immediately
leads to the following corollary.

\hfill

%%%%%%%%%%%%%%%%%%%%%%%%%%%%%%%%%%%%%%%%%%%%%%%%%%%%%%%%%%%%%%
\corollary\label{_vanishing_for_cone_Corollary_}
Let $(M,I)$ be a compact, irreducible,
holomorphically symplectic K\"ahler manifold, 
$\dim_{\Bbb C}M=2n$, and 
$L$ a holomorphic line bundle on $(M,I)$
with $c_1(L)\notin\bar{\cal K}\check{\;}$. Then
the holomorphic cohomology $H^i(M,L)$ are zero
for all $i<n$.

\endproof

\hfill

Now we can prove the main result of this paper.

\hfill

%%%%%%%%%%%%%%%%%%%%%%%%%%%%%%%%%%%%%%%%%%%%%%%%%%%%%%%%%%%%%%
\theorem\label{_main_vanishing_after_cones_Theorem_}
Let $(M,I,J,K)$ be a compact, irreducible hyperk\"ahler
manifold, and $L$ a holomorphic line bundle on $(M,I)$ with
$c_1(L)\neq 0$. Then one of the following holds.
\begin{description}
\item[(i)] $c_1(L)\in \bar{\cal K}\check{\;}$; then
$H^i(L) =0$ for $i> \frac{\dim_\C M}{2}$.
\item[(ii)] $c_1(L)\in -\bar{\cal K}\check{\;}$; then
$H^i(L) =0$ for $i< \frac{\dim_\C M}{2}$.
\item[(iii)] $c_1(L)$ does not lie in
$-\bar{\cal K}\check{\;}\cup \bar{\cal K}\check{\;}$; then
$H^i(L) =0$ for $i\neq \frac{\dim_\C M}{2}$.
\end{description}
{\bf Proof:} Denote $\frac{\dim_\C M}{2}$ by $n$.
\ref{_main_vanishing_after_cones_Theorem_}
(iii) is a direct consequence of
\ref{_vanishing_for_cone_Corollary_}.
Indeed, in this case 
\[ H^i(L) =0\text{\ \ for \ \ }i< n\]
and 
\begin{equation}\label{_coho_L^*_vanish_Equation_}
H^i(L^*) =0\text{\ \ for \ \ }i< n,
\end{equation}
because the Chern classes of both $L$ and $L^*$ 
do not lie in $\bar{\cal K}\check{\;}$.
However, by Serre's duality, 
\eqref{_coho_L^*_vanish_Equation_}
is equivalent to 
\[ H^i(L) =0\text{\ \ for \ \ }i> n.\]

Let us prove \ref{_main_vanishing_after_cones_Theorem_}
(i). Since 
$c_1(L)\in \bar{\cal  K}\check{\;}$,
we may assume that $q(c_1(L), \omega)\geq 0$
for all $\omega\in {\cal K}$. 
Unless $q(c_1(L), \omega)=0$ for
all K\"ahler classes $\omega$,
the assertion of  \ref{_main_vanishing_after_cones_Theorem_}
(i) is obtained from 
\ref{_vanishing_for_q(L,omega)>0_Proposition_}. However,
if $q(c_1(L), \omega)=0$ for
all K\"ahler classes, $c_1(L)=0$, because
$q$ is non-degenerate and the K\"ahler classes
generate $H^{1,1}_I(M,\R)$. \ref{_main_vanishing_after_cones_Theorem_}
(i) is obtained from (ii) by Serre's duality.
\endproof

\hfill

The classes  $\eta\notin-\bar{\cal K}\check{\;}\cup \bar{\cal
  K}\check{\;}$ can be also characterized as follows.

\hfill

%%%%%%%%%%%%%%%%%%%%%%%%%%%%%%%%%%%%%%%%%%%%%%%
\claim\label{_primi_and_q_Claim_}
Let $(M,I,J,K)$ be a compact, irreducible hyperk\"ahler
manifold, and $\eta\in H^{1,1}_I(M,\R)$
a non-zero cohomology class. Then  the following
conditions are clearly equivalent.
\begin{description}
\item[(i)] $\eta\notin-\bar{\cal K}\check{\;}\cup \bar{\cal
  K}\check{\;}$
\item[(ii)] $q(\eta, \omega_1)>0$ and $q(\eta, \omega_2)<0$ 
for some K\"ahler forms $\omega_1$, $\omega_2$ on $(M,I)$
\item[(iii)] $\eta$ is primitive with respect
to some K\"ahler form on $(M,I)$; or, equivalently,
$q(\eta,\omega)=0$ (see \ref{_BBF_and_primitive_Claim_}).
\item[(iv)] The class $\eta$ is $SU(2)$-invariant with respect
to some hyperk\"ahler structure $(I, J', K')$ on $M$.
\end{description}
{\bf Proof:} The equivalence of (i) and (ii) is clear.
The equivalence of (iii) and (iv) is implied by
\ref{_BBF_and_primitive_Claim_}. The implication (ii) $\Rightarrow$ (iii) is
clear, because the K\"ahler cone is connected, hence
from $q(\eta, \omega_1)>0$ and $q(\eta, \omega_2)<0$ 
it follows that $q(\eta,\omega_3)=0$ for some K\"ahler
form. Finally, (iii) $\Rightarrow$ (ii) is obtained as
follows: given a K\"ahler class $\omega$, with
$q(\eta,\omega)=0$, take a neighbourhood $U$ of $\omega$
in the K\"ahler cone. The function $U\stackrel v \arrow \R$,
$v(\omega')=q(\eta, \omega')$ is non-zero and 
linear, hence it takes positive and negative 
values in any open neighbourhood of $\omega$.
\endproof

%%%%%%%%%%%%%%%%%%%%%%%%%%%%%%%%%%%%%%%%%%%%%%%%%%%%%%%%%%%%%%
\subsection{Cohomology vanishing for vector bundles of
arbitrary rank}
%%%%%%%%%%%%%%%%%%%%%%%%%%%%%%%%%%%%%%%%%%%%%%%%%%%%%%%%%%%%%%

A version of \ref{_main_vanishing_after_cones_Theorem_}
can be stated for holomorphic bundles of arbitrary rank,
as follows. 

\hfill

%%%%%%%%%%%%%%%%%%%%%%%%%%%%%%%%%%%%%%%%%%%%%%%%%%%%%%%%%%%%%%
\theorem\label{_vanishing_arb_rank_Theorem_}
Let $(M,I,J,K)$ be a compact, irreducible hyperk\"ahler
manifold, $L$ a holomorphic line bundle on $(M,I)$ with
$c_1(L)\neq 0$, and $B$ an arbitrary holomorphic vector
bundle on $(M,I)$. Then there exists a sufficiently
big number $N_0$, such that for any integer $N>N_0$
one of the following holds.
\begin{description}
\item[(i)] $c_1(L)\in \bar{\cal K}\check{\;}$; then
$H^i(L^N\otimes B) =0$ for $i> \frac{\dim_\C M}{2}$.
\item[(ii)] $c_1(L)\in -\bar{\cal K}\check{\;}$; then
$H^i(L^N\otimes B) =0$ for $i< \frac{\dim_\C M}{2}$.
\item[(iii)] $c_1(L)$ does not lie in
$-\bar{\cal K}\check{\;}\cup \bar{\cal K}\check{\;}$; then
$H^i(L^N\otimes B) =0$ for $i\neq \frac{\dim_\C M}{2}$.
\end{description}
{\bf Proof:} The proof of \ref{_vanishing_arb_rank_Theorem_}
is similar to the Kodaira-Nakano vanishing for vector
bundles of arbitrary rank. The same argument as used
in the proof of \ref{_main_vanishing_after_cones_Theorem_}
can be employed to deduce \ref{_vanishing_arb_rank_Theorem_}
from the following statement.

\hfill

%%%%%%%%%%%%%%%%%%%%%%%%%%%%%%%%%%%%%%%%%%%%%%%%%%%%%%%%%%%%%%
\proposition\label{_vanishing_arb_rank_C_1(L)_pos_Proposition_}
Let $(M,I,J,K)$ be a compact, irreducible hyperk\"ahler
manifold, $L$ a holomorphic line bundle on $(M,I)$ with
$c_1(L)\neq 0$, and $B$ an arbitrary holomorphic vector
bundle on $(M,I)$. Assume that $q(c_1(L),\omega)>0$,
where $\omega$ is the K\"ahler form of $(M,I)$.
Then there exists a sufficiently
big number $N_0$, such that for any integer $N>N_0$
 $H^i(L^N\otimes B) =0$ for all $i> \frac{\dim_\C M}{2}$.

\hfill

{\bf Proof:} To prove
\ref{_vanishing_arb_rank_C_1(L)_pos_Proposition_},
we use the formula \eqref{_Laplacians_K-N_ide_Equation_}
again:
\begin{equation}\label{_Laplacian_ide_on_B_Equation_}
\Delta_{\bar\6}-\Delta_{\bar\6_J}= - [ \Theta_+, \Lambda_{\bar\Omega}],
\end{equation}
where $\Delta_{\bar\6}$, $\Delta_{\bar\6_J}$ 
are the Laplacians on $L^N\otimes B$, and
$\Theta_+$ is the 
$\Lambda^{1,1}_{+,I}(M)\otimes \End(L^N\otimes B)$-part
of the curvature of $L^N\otimes B$, considered
as an operator on $\Lambda^{(0,*)}(M)\otimes (L^N\otimes B)$
as in the proof of the quaternionic Dolbeault
Kodaira-Nakano identity
(\ref{_Kodaira_Nakano_Laplacians_Theorem_}).
Since the curvature is additive on tensor product, we have
\[
\Theta_+ = \Theta_B+ N \Theta_L,
\]
where $\Theta_B$, $\Theta_L$ are
$\Lambda^{2}_{+}(M)$-parts
of the curvatures of $B$ and $L$. 
The same argument as used in the proof of
\eqref{_Theta_+_for_line_bu_kahle_Equation_} implies that
$\Theta_L= \lambda\bar\Omega$, where
$\lambda=\frac{q(c_1(L),\omega)}{q(\omega,\omega)}$.
Then, as the Lefschetz formula (\ref{_SL(2)_Proposition_}) implies,
\[ - [ \Theta_+, \Lambda_{\bar\Omega}] =- [ \Theta_B,
  \Lambda_{\bar \Omega}]
+ V,
\]
where $V$ is a scalar operator acting on $(0,i)$-forms
as $\lambda (i-n)N$, $n=\frac{\dim_C M}2$. Clearly,
$- [ \Theta_B,
  \Lambda_{\bar \Omega}] + V$ is positive definite 
for $N$ sufficiently big and $i>n$. From 
\ref{_positivity_Lapla_Remark_}
we obtain immediately that $\ker \Delta_{\bar\6}=0$
whenever $- [ \Theta_B,
  \Lambda_{\bar \Omega}] + V$ is positive definite. 
This proves
\ref{_vanishing_arb_rank_C_1(L)_pos_Proposition_}.
\ref{_vanishing_arb_rank_Theorem_} is proven.
\endproof

%%%%%%%%%%%%%%%%%%%%%%%%%%%%%%%%%%%%%%%%%%%%%%%%%%%%%%%%%%%%%%

\section{Vanishing of cohomology and nef-classes with
  $q(\eta,\eta)=0$}
\label{_vanishing_and_nef_Section_}

%%%%%%%%%%%%%%%%%%%%%%%%%%%%%%%%%%%%%%%%%%%%%%%%%%%%%%%%%%%%%%

%%%%%%%%%%%%%%%%%%%%%%%%%%%%%%%%%%%%%%%%%%%%%%%%%%%%%%%%%%%%%%
\subsection{Nef classes}
%%%%%%%%%%%%%%%%%%%%%%%%%%%%%%%%%%%%%%%%%%%%%%%%%%%%%%%%%%%%%%

The following immensely important theorem was proven by
J.-P. Demailly and M. Paun.

\hfill

%%%%%%%%%%%%%%%%%%%%%%%%%%%%%%%%%%%%%%%%%%%%%%%%%%%%%%%%%%%%%%
\theorem\label{_Demailly_Paun_Theorem_}
(\cite{_Demailly_Paun_})
Let $M$ be a compact K\"ahler manifold,
and ${\cal X}$ the set of all closed analytic subvarieties
$X\subset M$ of positive dimension.
Consider the set of all $(1,1)$-classes
\[
\tilde {\cal K} := \{ \eta\in H^{1,1}(M)\cal H^2(M, \R) \ \ | \ \
\forall X\in {\cal X}, \ \ \int_X\eta^{\dim X}>0 \}
\]
Then the K\"ahler cone of $M$ coincides with
one of the connected components of $\tilde {\cal K}$.
\endproof

\hfill

%%%%%%%%%%%%%%%%%%%%%%%%%%%%%%%%%%%%%%%%%%%%%%%%%%%%%%%%%%%%%%
\remark
The converse assertion is trivial: if $\eta$ is a K\"ahler
class, then $\int_X\eta^i>0$ for all analytic cycles 
$X\subset M$, $\dim X=i$.

\hfill

%%%%%%%%%%%%%%%%%%%%%%%%%%%%%%%%%%%%%%%%%%%%%%%%%%%%%%%%%%%%%%
\definition
Let $M$ be a K\"ahler manifold, and $\eta\in H^{1,1}(M)$ a 
real $(1,1)$-class. Then $\eta$ is
called {\bf nef} (numerically effective) if $\eta$
belongs to a closure $\bar {\cal K}$ of the K\"ahler cone
${\cal K}$ of $M$. The closure $\bar {\cal K}$ 
is called {\bf the nef cone}. A nef line bundle on $M$
is a line bundle with $c_1(L)$ nef; a nef divisor $D$
is one with nef cohomology class.

%%%%%%%%%%%%%%%%%%%%%%%%%%%%%%%%%%%%%%%%%%%%%%%%%%%%%%%%%%%%%%
\subsection{Nef classes on hyperk\"ahler manifolds}
%%%%%%%%%%%%%%%%%%%%%%%%%%%%%%%%%%%%%%%%%%%%%%%%%%%%%%%%%%%%%%

Consider a compact, irreducible 
hyperk\"ahler manifold $(M,I,J,K)$. 
Let $L$ be a holomorphic line bundle on $(M,I)$
which is nef and satisfies 
\[ q(c_1(L),c_1(L))=0.\] It was conjectured
(\cite{_Huybrechts_}, \cite{_Sawon_})
that $L$ is base point free, that is,
defines a holomorphic map 
\begin{equation}\label{_map_to_P^n_Equation_}
(M,I) \arrow {\Bbb P} H^0(L^N)
\end{equation}
for $N$ sufficiently big. If this is true, then 
\eqref{_map_to_P^n_Equation_} is a Lagrangian
fibration onto its image (\cite{_Matsushita:fibred_}).
A special case of this conjecture was recently proven by
D. Matsushita (\cite{_Matsushita:nef_}). This motivates
our interest in the geometry of nef-classes 
satisfying $q(\eta,\eta)=0$.

\hfill

%%%%%%%%%%%%%%%%%%%%%%%%%%%%%%%%%%%%%%%%%%%%%%%%%%%%%%%%%%%%%%
\proposition\label{_defo_nef_obtain_outside_cones_Proposition_}
Let $(M,I,J,K)$ be a compact, irreducible hyperk\"ahler manifold,
$\eta\in H^{1,1}_I(M, \R)$ a non-zero nef class on $(M,I)$,
satisfying $q(\eta,\eta)=0$, and $\omega$ a rational
K\"ahler class on $(M, I)$. Then
\begin{description}
\item[(i)] $q(\omega, \eta)>0$

\item[(ii)] Choose a positive real number 
$\epsilon< \frac{q(\eta,\omega)}{q(\omega,\omega)}$.
Then $\eta-\epsilon\omega$ lies outside of 
$\bar{\cal K}\check{\;}\cup -\bar{\cal K}\check{\;}$.
\end{description}
{\bf Proof:} 
\ref{_defo_nef_obtain_outside_cones_Proposition_} (i)
is clear. Indeed, if $q(\omega, \eta)=0$, then
$q(\eta, \eta)<0$, because the form $q$ has signature
$(+, -, -, ... -)$ on $H^{1,1}(M, \R)$. On the other
hand, $q(\omega, \eta)\geq 0$, because $\eta$ lies in
the closure of the K\"ahler cone.

Let us prove
  \ref{_defo_nef_obtain_outside_cones_Proposition_} (ii).
Since $\epsilon< \frac{q(\eta,\omega)}{q(\omega,\omega)}$,
the number
\begin{equation}\label{_lambda_via_q_Equation_}
\lambda:= q(\omega, \eta-\epsilon\omega)
\end{equation}
is positive. Therefore, 
$\eta-\epsilon\omega\notin -\bar{\cal K}\check{\;}$.
To prove \ref{_defo_nef_obtain_outside_cones_Proposition_}
(ii), it remains to find a K\"ahler class $\omega'$
which satisfies $q(\omega', \eta-\epsilon\omega)<0$.

For any $\delta>0$, $\eta+\delta\omega$ is a K\"ahler
class, as follows from \ref{_Demailly_Paun_Theorem_}.
Choose a positive number
$\delta<\frac{\lambda\epsilon}{q(\eta\omega)}$,
where $\lambda$ is the number defined in 
\eqref{_lambda_via_q_Equation_}. Then
\[
q(\eta+\delta\omega, \eta-\epsilon\omega) =
-\epsilon\lambda +\delta q(\omega, \eta) = 
q(\omega, \eta)\left(\delta-\frac{\lambda\epsilon}{q(\eta,\omega)}\right)<0.
\]
\ref{_defo_nef_obtain_outside_cones_Proposition_} (ii) is proven.
\endproof

%%%%%%%%%%%%%%%%%%%%%%%%%%%%%%%%%%%%%%%%%%%%%%%%%%%%%%%%%%%%%%

\subsection{A vanishing theorem and its applications}

%%%%%%%%%%%%%%%%%%%%%%%%%%%%%%%%%%%%%%%%%%%%%%%%%%%%%%%%%%%%%%

From \ref{_defo_nef_obtain_outside_cones_Proposition_},
the following theorem is apparent.

\hfill

%%%%%%%%%%%%%%%%%%%%%%%%%%%%%%%%%%%%%%%%%%%%%%%%%%%%%%%%%%%%%%
\theorem\label{_vanishi_for_neigh_of_nef_Theorem_}
Let $(M,I,J,K)$ be a compact, irreducible hyperk\"ahler manifold,
$\dim_{\Bbb H}M=n$,
and $L$ a non-trivial holomorphic bundle on $(M,I)$
which is nef and satisfies $q(c_1(L),c_1(L)) =0$.
Consider an ample line bundle $H$ on $(M,I)$.
Then there exists $N_0$ such that for all 
integers $N>N_0$,
\begin{equation}\label{_vanishi_for_nef-minus-maple_Equation_}
H^i(L^N\otimes H^*)=0 \ \ \text{for} \ \  i\neq n.
\end{equation}
{\bf Proof:} Let $N_0= \frac 1 \epsilon$, where 
\[ \epsilon=\frac{q(c_1(L), c_1(H))}{q(c_1(H), c_1(H))}.
\]
Then $Nc_1(L)-c_1(H) \notin 
\bar{\cal K}\check{\;}\cup -\bar{\cal K}\check{\;}$
as follows from
\ref{_defo_nef_obtain_outside_cones_Proposition_}.
The vanishing of
\eqref{_vanishi_for_nef-minus-maple_Equation_}
then follows from
\ref{_main_vanishing_after_cones_Theorem_}.
\endproof

\hfill

\ref{_vanishi_for_neigh_of_nef_Theorem_} has an immediate
corollary.

\hfill

%%%%%%%%%%%%%%%%%%%%%%%%%%%%%%%%%%%%%%%%%%%%%%%%%%%%%%%%%%%%%%
\corollary \label{_restri_to_divi_Corollary_}
Let $(M,I,J,K)$ be a compact, irreducible hyperk\"ahler
manifold, $\dim_{\Bbb H}(M)>1$, 
$L$ a non-trivial holomorphic bundle on $(M,I)$
which is nef and satisfies $q(c_1(L),c_1(L)) =0$,
and $D$ an ample divisor on $(M,I)$. Then, for
sufficiently big $N> N_0$, the natural
restriction map \[ H^0(L^N) \arrow H^0(L^N\restrict D)\]
is surjective.

\hfill

{\bf Proof:}
The following exact sequence is well known
\[
0 \arrow L^N(-D) \arrow L^N \arrow L^N\restrict D\arrow 0.
\]
By \ref{_main_vanishing_after_cones_Theorem_},
$H^1(L^N(-D))=0$. Then the long exact sequence of
cohomology gives
\[
0 \arrow H^0(L^N(-D)) \arrow H^0(L^N) \arrow 
  H^0\left(L^N\restrict D\right)\arrow 0
\]
This proves \ref{_restri_to_divi_Corollary_}. \endproof

\hfill

\ref{_restri_to_divi_Corollary_} can be generalized as
follows.

\hfill

%%%%%%%%%%%%%%%%%%%%%%%%%%%%%%%%%%%%%%%%%%%%%%%%%%%%%%%%%%%%%%
\theorem\label{_comple_interse_surj_Theorem_}
Let $(M, I, J, K)$ be a irreducible hyperk\"ahler manifold,
and $X\subset (M,I)$ a subvariety of dimension
$\dim_\C X>\frac 1 2 \dim_\C M$. Assume that
$X$ is a complete intersection of  ample divisors.
Consider a holomorphic line bundle $L$ on
$(M,I)$ with $c_1(L)$ nef and
$q(c_1(L), c_1(L))=0$. Then the natural
restriction map is surjective on holomorphic sections:
\[
H^0(L^N) \arrow H^0\left(L^N\restrict X\right)\arrow 0,
\]
for a sufficiently big power of $L$.

\hfill

{\bf Proof:} Let $X = \bigcup_{i=1}^k H_i$, where
$k=\codim X$, and $H_i$ are ample divisors. Consider the
Koszul resolution of $L^N\restrict X$,
\begin{equation}\label{_Koszul_Equation_}
\begin{aligned}
0 \arrow & L^N(-H_1-H_2-\dots -H_k) \arrow \dots \\
\arrow & \bigoplus_{i>j} L^N(-H_i-H_j)\arrow 
\bigoplus_i L^N(-H_i)\\
\arrow & L^N \arrow L^N\restrict
X\arrow 0.
\end{aligned}
\end{equation}
By \ref{_vanishi_for_neigh_of_nef_Theorem_},
the cohomology of the terms $L^N(-H_i-H_j- ...)$
of \eqref{_Koszul_Equation_} vanish up to degree $n$.
Therefore, the $E_0^{p,q}$-term of the associated
spectral sequence looks as follows
 \begin{equation}\label{_E^0_array_Equation_}{\tiny
\begin{array}{ccccc}
H^n(L^N(-H_1- ... -H_k))& \dots  & H^n\bigg(\bigoplus_i L^N(-H_i)\bigg) & H^n(L^N) & 
                        H^n\left(L^N\restrict X\right)
\\[4mm] 0 & \dots &  0 & H^{n-1}(L^N) & 
                        H^{n-1}\left(L^N\restrict X\right)
\\ \vdots & {\ }  &\vdots & \vdots & \vdots
\\ 0 & \dots &  0 & H^{0}(L^N) & 
                        H^{0}\left(L^N\restrict X\right)
\end{array}}
\end{equation}
From \eqref{_E^0_array_Equation_} it is clear that the
only non-trivial differential of
\eqref{_E^0_array_Equation_} mapping to
$H^{0}\left(L^N\restrict X\right)$ is 
\begin{equation}\label{_d_1_restri_Equation_}
   d_1:\; H^{0}(L^N)\arrow H^{0}\left(L^N\restrict X\right)
\end{equation}
which is identified with the restriction map. Since the
complex \eqref{_Koszul_Equation_} is exact,
the spectral sequence \eqref{_E^0_array_Equation_}
converges to zero. Therefore, the differential
\eqref{_d_1_restri_Equation_} is surjective.
This proves \ref{_comple_interse_surj_Theorem_}.
\endproof

\hfill

{\bf Acknowledgements:} I am grateful to E. Amerik and
D. Markushevich for insightful comments and observations,
and to the referee for useful suggestions and corrections.

%%%%%%%%%%%%%%%%%%%%%%%%%%%%%%%%%%%%%%%%%%%%%%%%%%%%%%%%%%%%%%
{\scriptsize

\hfill

\noindent {\sc Misha Verbitsky\\
University of Glasgow, Department of Mathematics, \\
15 University Gardens, Glasgow G12 8QW, Scotland.}\\
\   \\
{\sc  Institute of Theoretical and
Experimental Physics \\
B. Cheremushkinskaya, 25, Moscow, 117259, Russia }\\
\  \\
\tt verbit@maths.gla.ac.uk, \ \  verbit@mccme.ru 
}% end of small

\end{document}